% Last Edit 24 February 1997

%            This document is written in Plain TeX
%            The macros: prepictex.tex, pictex.tex, and postpictex.tex are also
%             required for the full compilation of the document.

\magnification=1100
\overfullrule0pt

\input prepictex
\input pictex
\input postpictex
\input amssym.def

% ********************* Definitions ************************************

%\def\widetilde{\mathaccent"0365 }
\def\qed{\hbox{\hskip 1pt\vrule width4pt height 6pt depth1.5pt \hskip 1pt}}
\def\mapleftright#1{\smash{
   \mathop{\longleftrightarrow}\limits^{#1}}}

\def\CC{{\Bbb C}}
\def\RR{{\Bbb R}}
\def\ZZ{{\Bbb Z}}
\def\cA{{\cal A}}
\def\cB{{\cal B}}
\def\cC{{\cal C}}
\def\cF{{\cal F}}

\def\gp{{\goth p}}

% ********************* FONTS ************************************

\font\smallcaps=cmcsc10
\font\titlefont=cmr10 scaled \magstep1

\font\sectionfont=cmbx10
\font\tinyrm=cmr10 at 8pt

% ******************** SECTION HEADERS ***************************

\newcount\sectno
\newcount\subsectno
\newcount\resultno

\def\section #1. #2\par{
\sectno=#1
\resultno=0
\bigskip\noindent{\sectionfont #1.  #2}~\medbreak}

\def\subsection #1\par{ \global\advance\resultno by 1
\bigskip\noindent{\bf (\the\sectno.\the\resultno)  #1} \ \ }

%******************* MATHEMATICAL LABELS **************************

\def\prop{\bigskip\noindent{\bf Proposition. }\sl}
\def\propno{ \global\advance\resultno by 1
\bigskip\noindent{\bf (\the\sectno.\the\resultno)  Proposition. }\sl}

\def\remark{\bigskip\noindent{\sl Remark. }}

\def\conj{\bigskip\noindent{\bf Conjecture. }\sl}
\def\exampleno{ \global\advance\resultno by 1
\bigskip\noindent{\bf (\the\sectno.\the\resultno)  Example. }\sl}
\def\conjno{ \global\advance\resultno by 1
\bigskip\noindent{\bf (\the\sectno.\the\resultno)  Conjecture. }\sl}
\def\thmno{ \global\advance\resultno by 1
\bigskip\noindent{\bf (\the\sectno.\the\resultno)  Theorem. }\sl}
\def\thm{\bigskip\noindent{\bf Theorem. }\sl}
\def\endthm{\rm\bigskip}

\def\endexample{\rm\bigskip}
\def\endprop{\rm\bigskip}

\def\pf{\rm\bigskip\noindent{\it Proof. }}
\def\endpf{\qed\hfil\bigskip}

%*************** EQUATIONS WITH NUMBERS **************

\def\formula{\global\advance\resultno by 1
\eqno{(\the\sectno.\the\resultno)}}
\def\formulano{\global\advance\resultno by 1 (\the\sectno.\the\resultno)}
\def\tableno{\global\advance\resultno by 1
\the\sectno.\the\resultno. }
\def\lformula{\global\advance\resultno by 1
\leqno(\the\sectno.\the\resultno)}

%********** DATING ******************************************
\def\monthname {\ifcase\month\or January\or February\or March\or April\or
May\or June\or
July\or August\or September\or October\or November\or December\fi}

\newcount\mins  \newcount\hours  \hours=\time \mins=\time
\def\now{\divide\hours by60 \multiply\hours by60 \advance\mins by-\hours
     \divide\hours by60         % NOTE: \divide only gives integer answers.
     \ifnum\hours>12 \advance\hours by-12
       \number\hours:\ifnum\mins<10 0\fi\number\mins\ P.M.\else
       \number\hours:\ifnum\mins<10 0\fi\number\mins\ A.M.\fi}

%**************** PAGE HEADERS *************************

\nopagenumbers
\def\runningtitle{\smallcaps standard young tableaux for root systems}
\headline={\ifnum\pageno>1\eoheadline\else\firstheadline\fi}
\def\names{\smallcaps arun ram}
\def\firstheadline{\noindent Preliminary Draft \hfill  \today}
\def\firstheadline{}
\def\eoheadline{\ifodd\pageno\oddheadline\else\evenheadline\fi}
\def\oddheadline{\tenrm\hfil\runningtitle\hfil\folio}
\def\evenheadline{\tenrm \folio\hfil{\names}\hfil}

%**************** TITLE *************************
\vphantom{$ $}  %My kludge to get the first page to move down a bit
\vskip.75truein

\centerline{\titlefont Standard Young tableaux for finite root systems}
%\smallbreak
%\centerline{\titlefont IWAHORI-HECKE ALGEBRA OF TYPE A${}_{\bf n-1}$}
\bigskip
\centerline{\rm Arun Ram${}^\ast$}
\centerline{Department of Mathematics}
\centerline{Princeton University}
\centerline{Princeton, NJ 08544}
\centerline{{\tt rama@math.princeton.edu}}
\centerline{Preprint: March 26, 1998}

\footnote{}{\tinyrm ${}^\ast$ Research supported in part by National
Science Foundation grant DMS-9622985 and a Postdoctoral Fellowship
at the Mathematical Sciences Research Institute.}

\bigskip

%**************** ABSTRACT *************************
\noindent{\bf Abstract.}
The study of representations of affine Hecke algebras has led to a new
notion of shapes and standard Young tableaux which works for the root system of
any finite Coxeter group.  This paper is completely independent of
affine Hecke algebra theory and is purely combinatorial.
We define generalized shapes and standard Young
tableaux and show that these new objects coincide with the classical ones
for root systems of Type A.  The classical notions of
conjugation of shapes, ribbon shapes, axial distances,
and the row reading and column reading standard tableaux,
have natural generalizations to the root system case.
In the final section we give an interpretation of the shapes and
standard tableaux for root systems of Type C
which is in a form similar to classical theory of shapes and standard tableaux.

\section 0. Introduction

In my recent work on representations of affine Hecke algebras
[Ra1] I have been led to a generalization of standard Young
tableaux. These generalized tableaux are important in the context of
representation theory because the standard tableaux model the internal
structure of irreducible representations of the affine Hecke algebra.
In fact, most of the time the number of tableaux of a given 
shape is the same as the dimension of the corresponding irreducible 
representation of the affine Hecke algebra.

In this paper I introduce and study generalized shapes and standard
tableaux purely combinatorially.  The main theorem is that the
generalized standard tableaux of a given shape
describe the connected components of a certain graph, the calibration graph.  
It is this graph which is intimately connected to the structure of 
representations of affine Hecke algebras.

In the Type A case the generalized shapes can be
converted into ``placed configurations of boxes''.  
This conversion is nontrivial and is the subject of Section 3.
In the cases where this placed configuration of boxes is a placed skew shape
the generalized standard tableaux coincide with the classical standard
tableaux of a skew shape. The generalized skew shapes 
play a major role in the results on representations of affine Hecke 
algebras which are obtained in [Ra1]. 

In Section 1 I give definitions of
\smallskip
\itemitem{(a)}  skew shapes,
\itemitem{(b)}  ribbon shapes,
\itemitem{(c)}  axial distances,
\itemitem{(d)}  conjugation of shapes, and
\itemitem{(e)}  row reading and column reading tableaux,
\smallskip\noindent
in the generalized setting.  In Section 4 it is shown that 
these definitions yield the classical versions of these objects in the Type A
case. The last section of this paper explains how one can convert the
generalized shapes and standard tableaux for the Type C case
into configurations of boxes and fillings.
In this form the shapes and standard tableaux for Type C look similar
to the classical standard Young tableaux.

It is my hope that others will also take up the study of the
generalized shapes and standard tableaux introduced in this paper.
There are many more problems than there is time for solving them and
every combinatorial fact
which can be proved about these objects says something about the structure of
representations of affine Hecke algebras.
One hopes that everything that is known about classical standard Young tableaux
will have an analogue in this more general setting.  Although I have
uncovered some of these generalizations, there are many facets of classical 
tableaux theory which still need to be generalized. 

From a representation theoretic point of view, one expects that there might
exist generalizations of
\itemitem{(a)}  the Robinson-Schensted-Knuth correspondence,
\itemitem{(b)}  the Littlewood-Richardson coefficients,
\itemitem{(c)}  the Kostka-Foulkes polynomials,
\itemitem{(d)}  major index and descents of tableaux,
\itemitem{(e)}  charge of tableaux,
\itemitem{(f)}  Jacobi-Trudi formulas.
\smallskip\noindent
Any solutions to these problems would be extremely helpful for 
understanding the underlying representation theory.
It is possible that some of the generalizations
might be obtained simply by understanding how to do them 
for the Type C tableaux given
in the last section of this paper.  This approach is attractive since
the form of the Type C tableaux given in the last section looks so
similar to classical tableau theory.

\bigskip\noindent
{\it Remarks on the results in this paper}

\medskip
\item{(1)}
Recently, it has become clear to me that one of the reasons that skew Schur
functions play such a important role in the classical theory of symmetric
functions is because the skew shapes describe particularly well behaved
{\it irreducible} representations of the affine Hecke algebra of type A.  Since
these nice representations exist in all types it seems reasonable that the many
wonderful identities involving skew Schur functions should have general type
analogues in terms of the generalized skew shapes defined in this paper. 
An example  of a skew Schur function identity that has a particularly nice
generalization to all types is the identity in [Mac I \S5 Ex.~21b].

\item{(2)} A remark similar to (1) can be made concerning ribbon shapes.  
This special class
of shapes has a good generalization to all types and the representation theory
associated to ribbon shapes has special features [Ro], [Mat, 4.3.5].
This fact seems to give some philosophical ``reason'' why there is such 
an amazing theory of ribbon Schur functions. The theory of ribbon 
Schur functions has been
developed in the last decade by Lascoux, Leclerc, Thibon, Krob, Reutenauer,
Malvenuto and others [GK] [La1-2].
I am sure that there is much to learn about the representation theory
of affine Hecke algebras from what is already known in the ribbon (and
noncommutative) Schur function theory.  

\item{(3)}  The theory of generalized shapes gives rise to some strange new
shapes even in the type A case, see Section 3.6.  To my knowledge
these shapes have not been studied before but they do retain many of the
combinatorial properties that skew shapes have.  In particular,
standard tableaux make perfectly good sense for these shapes and 
these strange standard tableaux do have representation theoretic meaning.

\item{(4)}  There seem to be strong connections between the combinatorics 
in this paper, the combinatorics of the Shi arrangements (see [St1-2],
[AL], [He], [ST]) and the combinatorics
of sign types developed by Shi [Sh3].  These connections need to be better
understood.

\item{(5)}  I do {\it not} think that there is a connection between the 
generalized standard Young tableaux introduced in this paper and the 
generalized tableaux of P. Littelmann [Li1-2].  Littelmann's
analogue of tableaux are really a generalization of column strict
tableaux {\it not} of standard tableaux.  Column strict tableaux
give information about the representations of $GL_n(\CC)$
and standard tableaux give information about the representations
of the symmetric group $S_n$.  There is an analogous dichotomy
in the generalized case; Littelmann's generalized 
column strict tableaux model the representations of complex semisimple 
Lie groups and my generalized standard tableaux model the representations
of affine Hecke algebras.  At the moment, I do not believe
that there is any connection  between the representation theories of 
the complex Lie groups and the affine Hecke algebras 
which would allow one to transfer information from one side to the other
(except in the type $A$ case, where one has a Schur-Weyl type duality.) 

\bigskip\noindent
{\it Acknowledgements} 
 
\medskip
This paper is only a part of a large project [Ra1-3] [RR1-2] on
representations of affine Hecke algebras which I have been working on intensely
for about a year.  During that time I have benefited from conversations with
many people. To  choose only a few, there were discussions with
S. Fomin, M. Vazirani, L. Solomon, F. Knop and N. Wallach which played an
important role in my progress.  There were several times when
I tapped into J. Stembridge's fountain of useful knowledge about root systems. 
G. Benkart was a very patient listener on many occasions.  
R. Simion, T. Halverson, H. Barcelo, P. Deligne, and 
R. Macpherson all gave large amounts of time to let me tell them
my story and every one of these sessions was helpful to me in solidifying
my understanding.  I thank C. Kriloff for her amazing proofreading.

I single out Jacqui Ramagge with special thanks for everything she
has done to help with this project: from the most mundane typing and picture
drawing to deep intense mathematical conversations which helped to sort out
many pieces of this theory.  Her immense contribution is evident in
that some of the papers in this series on representations of affine Hecke
algebras are joint papers.  

A portion of this research was done during a semester stay at Mathematical
Sciences Research Institute where I was supported by a Postdoctoral
Fellowship.  I thank MSRI and National Science Foundation for support of
my research.

\section 1. Generalized shapes and standard Young tableaux

\subsection{Notations.}

Let $W$ be a finite Coxeter group and let $R$ be the root system of $W$.  
The root system $R$ spans a real vector space which we shall denote
$\RR^n$.  There is an inner product on $\RR^n$ via which $W$ is a 
group generated by reflections.  
%The roots in $R$ describe a (central) hyperplane arrangement
%$$\cA = \{ H_\alpha \ |\ \alpha\in R\},
%\qquad\hbox{where}\qquad
%H_\alpha = \{ x\in \RR^n \ |\ \langle x,\alpha\rangle = 0\}.$$
Fix a system $R^+$ of {\it positive roots} in $R$ and write 
$\alpha>0$ if $\alpha\in R^+$.   
If $w\in W$ the {\it inversion set} of $w$ is
$$R(w) = \{ \alpha> 0 \ |\ w\alpha<0 \,\}.$$
Let $\{\alpha_1,\ldots,\alpha_n\}$ be the {\it simple roots} in $R^+$ and let
$s_1,\ldots,s_n$ denote the corresponding {\it simple reflections} in $W$. 
The positive roots determine a {\it fundamental chamber}
$C=\{ x\in \RR^n\ |\ \langle x,\alpha\rangle >0
\hbox{\ for all $\alpha\in R^+$} \}.$ 
An element $x\in \RR^n$ is {\it dominant} if $x\in \overline C$, 
the closure of the chamber $C$.  
%An element $x\in \RR^n$ is 
%{\it regular} if $\langle x, \alpha\rangle \ne 0$ for all $\alpha\in R$,
%and is {\it integral} if $\langle x,\alpha\rangle\in \ZZ$ for all
%$\alpha\in R$. 

\subsection{Placed shapes and standard tableaux.}

If $\gamma\in \overline{C}$ define
$$
Z(\gamma)=\left\{ \alpha\in R^+ \ |\  \langle\gamma,\alpha\rangle=0
\right\}
\qquad\hbox{and}\qquad
P(\gamma)=\left\{ \alpha\in R \ |\  \langle\gamma,\alpha\rangle= 1
\right\}.
$$
The set $Z(\gamma)\cup -Z(\gamma)$ is a (parabolic) root subsystem of $R$ and is
generated by the simple roots that it contains.  The stabilizer $W_\gamma$ 
of $\gamma$ in $W$ is the Weyl group of the subsystem $Z(\gamma)$ and
the set $P(\gamma)$ is stable under the action of $W_\gamma$.

\remark
The sets $Z(\gamma)$ and $P(\gamma)$ appear in  
Heckman and Opdam [HO, Definition 2.1] and implicitly in other works
[K1, Theorem 2.2].
\endexample

A {\it placed shape} is a pair $(\gamma,J)$ where $\gamma\in \overline{C}$
and $J\subseteq P(\gamma)$.  
A {\it standard tableau of shape $(\gamma,J)$} is 
an element $w\in W$ such that
$$R(w)\cap Z=\emptyset \qquad\hbox{and}\qquad R(w)\cap P=J.$$
In Theorem 3.5 we shall show that 
this is a generalization of the classical notion of standard Young tableau. 
Let 
$$\cF^{(\gamma,J)}=\{\hbox{ standard tableaux of shape $(\gamma,J)$} \}.$$
%If $\gamma\in \overline{C}$ and $\mu\in W\gamma$ then there is a 
%unique element $w\in W$, of minimal length, such that $w\gamma=\mu$.
%The placed shape determined by $\gamma$ and $\mu$ is
%$$\gamma/\mu= (\gamma,J(\mu)),
%\qquad\hbox{where $J(\mu) = R(w)\cap P(\gamma)$.}
%$$

\medskip\noindent
{\sl Examples}
\item{(1)}
If $\gamma$ is a generic element of $\RR^n$ then 
$Z(\gamma)=P(\gamma)=\emptyset$.
In this case the only possibility for $J$ is $J=\emptyset$ and
$\cF^{(\gamma,\emptyset)}=W$.
\item{(2)}
Let $\rho={1\over2}\sum_{\alpha>0}\alpha^\vee$, where $\alpha^\vee =
2\alpha/(\alpha,\alpha)$. Then $Z(\rho)=\emptyset$ and
$P(\rho)=\{\alpha_1,\ldots,\alpha_n\}$. If $J\subseteq
\{\alpha_1,\ldots,\alpha_n\}$ we have 
$$\cF^{(\rho,J)}~=~\{ w\in W\ |\ D(w)=J\},
\qquad\hbox{where\quad  $D(w) = \{ \alpha_i \ |\ ws_i<w\}$}$$ 
is the {\it right descent set} of $w\in W$.

\subsection{A nonemptiness condition.}

I feel that the following conjecture should have a simple slick proof.  
The $\Longrightarrow$ direction is easy and it may be possible
to prove the other direction by using [Bou, Ch. VI \S 1, Prop. 22].
For Type A the conjecture is true and is a consequence of Theorem 4.5
of section 4.

\conj  Let $(\gamma,J)$ be a shape.  The set $\cF^{(\gamma,J)}$
is non-empty if and only if $J$ satisfies the condition
$$
\hbox{ {\sl  If $\beta\in J$, $\alpha\in Z$ and $\beta-\alpha\in R^+$
then $\beta-\alpha\in J$.}}
$$
\rm

\subsection{Placed skew shapes.}

If $\gamma\in \RR^n$ view $\gamma$ as the function on the root system
$R$ given by
$$\matrix{
\gamma\,\colon &R &\longrightarrow &\RR \cr
&\alpha &\longmapsto &\langle \gamma,\alpha\rangle \cr
}
$$
The element $\gamma$ is {\it regular} if 
$\langle \gamma,\alpha\rangle\ne 0$ for all $\alpha\in R$, and is {\it
integral} if $\langle \gamma,\alpha\rangle\in \ZZ$ for all $\alpha\in R$.
For each subset $K\subseteq\{\alpha_1,\ldots,\alpha_n\}$ let
$R_K$ be the root system generated by $K$ and let $W_K$ be
the Weyl group of $R_K$.
Let $\gamma\big\vert_K$ denote the function $\gamma$ restricted to $R_K$.

A placed shape $(\gamma,J)$ is a placed {\it skew shape} if
for all $w\in \cF^{(\gamma,J)}$,
\smallskip
\itemitem{(a)}  for each simple root $\alpha_i$, 
$w\gamma\big\vert_{\{\alpha_i\}}$ is regular, and
\smallskip
\itemitem{(b)}  for each pair of simple roots $\{\alpha_i,\alpha_j\}$ 
\itemitem{}
$\displaystyle{
\qquad\matrix{
\hbox{either}\qquad &\hbox{$w\gamma\big\vert_{\{\alpha_i,\alpha_j\}}$
is regular}\hfill \cr
\hbox{or}
&\hbox{$w\gamma\big\vert_{\{\alpha_i,\alpha_j\}}$ is in the
$W_{\{\alpha_i,\alpha_j\}}$-orbit of} \hfill \cr
&\hbox{the function $\tau$ given by
$\langle \tau,\alpha_i\rangle =1$ and
$\langle \tau,\alpha_j\rangle =0$,} \hfill \cr
&\hbox{where $\alpha_i$ is long and $\alpha_j$ is short.} \hfill \cr
} 
}
$

\remark
In the type A case, a definition of skew shape similar
to the one given above 
has also been given by Fomin [Fo] in connection with his
approach to the theory of seminormal representations of the symmetric group.
\endexample

\subsection{Ribbon shapes.}

A placed shape $(\gamma,J)$ is a {\it placed ribbon shape} if
$Z(\gamma)=\emptyset$. All placed ribbon shapes are placed skew shapes.

\smallskip\noindent
{\sl Example.}
Suppose that $(\gamma,J)$, is a ribbon shape and $\gamma$ is
integral. Then
$Z(\gamma)=\emptyset$, $P(\gamma)\subseteq \{\alpha_1,\ldots,\alpha_n\}$ and
$$\cF^{(\gamma,J)}~=~\{ w\in W\ |\ L(w)\cap P=J\},$$
where, as in Example 2 of (1.2), $L(w)$ is the left descent set of $w$.

\subsection{}

Although we shall not define the affine Hecke algebra
or discuss its representation theory in this paper it is important to note that
the definition of placed skew shape is motivated by the following theorem 
of [Ra1].

\thm  Let $R$ be the root system of a finite Weyl group and
let $\tilde H$ be the corresponding affine Hecke algebra.
There is a one-to-one correspondence between placed skew shapes $(\gamma,J)$
and irreducible really calibrated representations $\tilde H^{(\gamma,J)}$ 
of the affine Hecke algebra $\tilde H$.  Under this correspondence
$$\dim(\tilde H^{(\gamma,J)})=
(\hbox{ \# of standard tableaux of shape $(\gamma,J)$}).$$
\endthm

\subsection{Conjugation.}

Let $(\gamma,J)$ be a placed shape and let $W_\gamma$ be the stabilizer
of $\gamma$ in $W$.  Define the {\it conjugate} placed shape to be
$$
(\gamma,J)'=(-u\gamma,-u(P(\gamma)\setminus J)),
$$ 
where $u$ is the minimal length coset representative of
$w_0W_\gamma\in W/W_\gamma$ and $w_0$ is the longest element of $W$. 

It will be useful to note the following:
$$
\hbox{(a)\quad $P(-u\gamma)=-uP(\gamma)$,}
\qquad\qquad \hbox{(b)\quad $Z(-u\gamma)=uZ(\gamma)$,}
\qquad\qquad \hbox{(c)\quad $R(u) = R^+\backslash Z(\gamma)$.}
$$
The proofs are as follows:

\item{(a)}  Since $\gamma$ is dominant, $-u\gamma=-w_0\gamma$ is
dominant and thus $\langle -u\gamma, -u\alpha\rangle=1$
only if $-u\alpha>0$.
With this in mind $P(-u\gamma)=-uP(\gamma)$ follows from the equation
$\langle -u\gamma,-u\alpha\rangle =1 
\Longleftrightarrow \langle \gamma,\alpha\rangle = 1$.

\item{(b)}  Let $v\in W_\gamma$ such that $w_0=uv$.  
(By [Bou IV \S 1 Ex. 3], $v$ is unique.)
Then $R^+\supseteq -w_0Z(\gamma)=-uvZ(\gamma)=uZ(\gamma)$, and it follows 
that
$$Z(-u\gamma) 
=R^+ \cap \{ \alpha\in R\ |\ \langle u\gamma,\alpha\rangle=0\}
=R^+\cap (uZ(\gamma)\cup -uZ(\gamma))=uZ(\gamma).$$

\item{(c)}  Let $R^-=-R^+$ be the set of negative roots in $R$.
Let $v\in W_\gamma$ such that $w_0=uv$.  Then
$v$ is the longest element of $W_\gamma$ and $R(v)=Z(\gamma)$. Thus
$$
\eqalign{
R(u)&=
\{ \alpha\in R \ |\ \alpha\in R^+, w_0v\alpha\in R^-\} \cr
&=\{ \alpha\in R \ |\ \alpha\in R^+, v\alpha\in R^+\},
\quad\hbox{since $w_0R^-=R^+$,} \cr
&= R^+\backslash R(v) = R^+\backslash Z(\gamma). \cr
}
$$

\propno  Conjugation $(\gamma,J)\longleftrightarrow (\gamma,J)'$
is a well defined involution on placed shapes. 
\pf
\item{(a)}
The weight $-u\gamma = -uv\gamma = -w_0\gamma$ is dominant
(i.e. $-u\gamma\in \overline{C}$)
and  $-u(P(\gamma)\backslash J)\subseteq P(-u\gamma)$ since
$-uP(\gamma)=P(-u\gamma)$. This shows that $(\gamma,J)'$ is
well defined.
\item{(b)}
Write $w_0=uv$ where 
where $v$ is the longest element of $W_\gamma$.  Similarly,
write $w_0=u'v'$ where $u'$ is the minimal
length coset representative of $w_0W_{w_0\gamma}$ and $v'$ is
the longest element in $W_{w_0\gamma}$. 
Conjugation by $w_0$ is an involution on $W$ which takes simple reflections to
simple reflections and $W_{w_0\gamma}=w_0W_\gamma w_0$.  It follows that
$v'=w_0v w_0$.  This gives
$$u'u=(w_0v')(w_0v)=w_0w_0vw_0w_0v =1.$$
Then, using fact (1.7a) above,
$$
-u'(P(-u\gamma)\backslash(-u(P(\gamma)\backslash J)))
=-u'(-uP(\gamma)\backslash(-u(P(\gamma)\backslash J)))
=P(\gamma)\backslash(P(\gamma)\backslash J)=J.
\qquad\qquad\hbox{\qed}$$
\medskip

\noindent
In (4.4) we shall show that this
involution is a generalization of the classical conjugation 
operation on partitions.

\remark
In type $A$, the conjugation involution seems to coincide
with the duality operation for representations of $\gp$-adic $GL(n)$ defined by
Zelevinsky [Ze].  Zelevinsky's involution has been been studied
further in [MW] and [KZ] extended to general Lie type by Kato [K2] 
and Aubert [Au].

\propno
The conjugation of shapes involution extends to an involution on
standard tableaux given by
$$\matrix{
\cF^{(\gamma,J)} &\mapleftright{1-1} &\cF^{(\gamma,J)'} \cr
w &\longleftrightarrow &wu^{-1}\cr
}
$$
where $u$ is the minimal length coset representative of
the coset $w_0W_\gamma$.
\pf
Let $w\in \cF^{(\gamma,J)}$ and let $w'=wu^{-1}$.  We must show that
$R(w')\cap Z(-u\gamma)=\emptyset$ and 
$R(w')\cap P(-u\gamma)=-u(P(\gamma\setminus J))$.
Let $R^-=-R^+$ be the set of negative roots.
We shall use the facts in (1.7) freely.

\smallskip\noindent
(a)  
Since $R(w)\cap Z(\gamma)=\emptyset$, we get
$$
\eqalign{
u^{-1}R(wu^{-1})\cap Z(\gamma)
&= \{\beta\in R\ |\ u\beta\in R(wu^{-1}), \beta\in Z(\gamma) \} \cr
&=\{\beta\in R\ |\ u\beta\in R^+,wu^{-1}u\beta\in  R^-,\beta\in Z(\gamma) \}\cr
&= \{\beta\in R\ |\ \beta\in u^{-1}R^+,w\beta\in  R^-,\beta\in Z(\gamma) \}\cr
&= \{\beta\in R\ |\ \beta\in u^{-1}R^+,\beta\in R(w), \beta\in Z(\gamma) \}, 
\quad\hbox{since $Z(\gamma)\subseteq R^+$,} \cr
&= \{\beta\in R\ |\ \beta\in u^{-1}R^+, \beta\in R(w)\cap Z(\gamma) \}\cr
&=\emptyset,\cr
}
$$
and it follows that
$$R(w')\cap Z(-u\gamma)
=R(wu^{-1})\cap uZ(\gamma)
=u\left(u^{-1}R(wu^{-1})\cap Z(\gamma)\right)
=\emptyset.$$
(b)  Assume that $w\in W$ is such that $R(w)\cap P(\gamma)=J$.  
Then
$$
\eqalign{
-u^{-1}R(wu^{-1})\cap P(\gamma)
&= \{\beta\in R\ |\ -u\beta\in R(wu^{-1}), \beta\in P(\gamma) \} \cr
&= \{\beta\in R\ |\ -u\beta\in R^+, -wu^{-1}u\beta\in  R^-,
\beta\in P(\gamma) \} \cr
&= \{\beta\in R\ |\ u\beta\in R^-, w\beta\in  R^+, \beta\in P(\gamma) \} \cr
&= \{\beta\in R\ |\ \beta\in R(u), \beta\in  R^+\backslash R(w),
\beta\in P(\gamma) \}, \quad\hbox{since $P(\gamma)\subseteq R^+$} \cr
&= \{\beta\in R\ |\ \beta\in R^+\backslash Z(\gamma),
\beta\in  R^+\backslash R(w), \beta\in P(\gamma) \} \cr
&=\{\beta\in R\ |\ \beta\in R^+\backslash Z(\gamma),
\beta\in  P(\gamma)\backslash J \},
\quad\hbox{since $R(w)\cap P(\gamma)=J$,} \cr
&= P(\gamma)\backslash J,
\quad\hbox{since $Z(\gamma)$ and $P(\gamma)$ are disjoint.} \cr
}
$$
This yields
$$R(w')\cap P(-u\gamma)
=R(wu^{-1})\cap -u P(\gamma)
=-u\left( -u^{-1}R(wu^{-1})\cap P(\gamma) \right)
=-u\left(P(\gamma)\backslash J\right).\qquad\qquad\hbox{\qed}$$
\medskip

\subsection{Axial distances.}

Let $(\gamma,J)$ be a placed shape and let $w\in \cF^{(\gamma,J)}$
be a standard tableaux of shape $(\gamma,J)$.  Let
$\alpha\in R$ be a root.  The {\it $\alpha$-axial distance} 
for $w$ is the value
$$d_\alpha(w) = \langle w\gamma, \alpha\rangle.$$
These numbers are crucial to the construction of irreducible
representations of the affine Hecke algebra which is given in
[Ra1].  In (4.1) we shall see that they are analogues
of the axial distances used by A. Young [Y] in his constructions
of the irreducible representations of the symmetric group.

\subsection{Row reading and column reading tableaux.}

If $w\in W$ let $R(w)$ be the inversion set of $w$.
The {\it weak Bruhat order} is the partial order on $W$ given by
$$\hbox{$v\le w$\qquad if\qquad $R(v)\subseteq R(w)$.}
$$
This definition is not the usual definition of the weak Bruhat order but is
equivalent to the usual one by [Bj, Prop. 2].  

Let $(\gamma,J)$ be a placed shape.
A {\it column reading tableau} of shape $(\gamma,J)$ is a
minimal element of $\cF^{(\gamma,J)}$ in the weak Bruhat order.
A {\it row reading tableau} of shape $(\gamma,J)$ is a 
maximal element of $\cF^{(\gamma,J)}$ in the weak Bruhat
order.

A set of positive roots $K$ is {\it closed} if  
$\alpha,\beta\in K$, $\alpha+\beta\in R^+$ implies
that $\alpha+\beta\in K$.  The {\it closure} $\overline{K}$ of a subset
$K\subseteq R^+$ is the smallest closed subset of $R^+$ containing $K$.

%The following result is stated without proof in [Bj, Prop. 3].
%J. Stembridge has given an inductive proof which can be found in [Wa]. 
%
%\prop Let $K\subseteq R^+$.  Then there exists
%a $w\in W$ such that $R(w)=K$ if and only if
%$K$ is closed and $K^c$ is closed.
%\endprop

The following conjecture says that if $\gamma$ is
integral and $\cF^{(\gamma,J)}$ is nonempty then
$\cF^{(\gamma,J)}$ contains a unique column reading tableau and
a unique row reading tableau.  A proof of the conjecture in the
type A case is given in Theorem 4.5.

\conj  Let $(\gamma,J)$ be a placed shape such that $\gamma$ is 
dominant and integral.
If $\cF^{(\gamma,J)}\ne \emptyset$ then 
$$R(w_{\rm min})=\overline{J},
\qquad R(w_{\rm max})= \overline{(P(\gamma)\setminus J)\cup Z(\gamma)}^c,
\qquad\hbox{and}\qquad
\cF^{(\gamma,J)} = 
[w_{\rm min},w_{\rm max}],
$$
where 
$K^c$ denotes the complement
of $K$ in $R^+$
and $[w_{\rm min},w_{\rm max}]$ denotes the interval between $w_{\rm min}$ and
$w_{\rm max}$ in the weak Bruhat order. 
\endthm

%\remark  The following example shows that the assumption that
%the above conjecture is false without the integrality assumption
%on $\gamma$.  Let $R$ be the root system of type $C_2$ with simple roots
%$$\alpha_1=\varepsilon_1
%\qquad\hbox{and}\qquad
%\alpha_2=\varepsilon_2-\varepsilon_1,$$
%where $\{\varepsilon_1,\varepsilon_2\}$ is an orthonormal basis
%of $\RR^2$.  Let $\gamma\in \RR^2$ be defined by
%$$\langle \gamma, \alpha_1\rangle = 0,
%\qquad
%%\langle \gamma,\alpha_2\rangle = \hbox{${1\over 2}$}.$$
%Then $Z(\gamma)=\{\alpha_1\}$, $P(\gamma)=\{\alpha_1+2\alpha_2}$.
%If $J=P(\gamma)$ then the unique minimal element 
%$w_{\rm min}$ of $cF^{(\gamma,J)}$ has 
%$R(w_{\rm min})=\{\alpha_2,\alpha_1+2\alpha_2\}\ne \overline{J}=J$.
%\endexample

\section 2. Calibration graphs

The following graphs arise naturally in the study of representations
of affine Hecke algebras, see [Ro], [Rg, Prop. 3.5], and
[Ra1].
Let $\gamma\in \RR^n$.  The {\it calibration graph} $\Gamma(\gamma)$ 
is the graph with 
$$
\matrix{
\hbox{Vertices:} \quad  &W\gamma\hfil  \cr
\hbox{Edges:}  &w\gamma \longleftrightarrow s_iw\gamma,\quad
&\hbox{ if }\quad \langle w\gamma,\alpha_i\rangle \ne \pm1. \hfil \cr
}
$$

\thmno
Assume that $\gamma\in \RR^n$ is dominant and let $\Gamma(\gamma)$ be the 
corresponding calibration graph.
The connected components of $\Gamma(\gamma)$ are described by the
sets 
$$
\cF^{(\gamma,J)} \qquad\hbox{such that}\qquad J\subseteq P(\gamma)
\quad\hbox{and}\quad \cF^{(\gamma,J)}\neq\emptyset,
$$
where $\cF^{(\gamma,J)}$ is the set of standard tableaux of 
shape $(\gamma,J)$ defined in (1.2).
\endthm

\subsection{}

This theorem will become almost obvious once we change our point of
view. The root system $R$ determines a 
(central) hyperplane arrangement
$$\cA = \{ H_\alpha \ |\ \alpha\in R\},
\qquad\hbox{where}\qquad
H_\alpha = \{ x\in \RR^n \ |\ \langle x,\alpha\rangle = 0\}.$$
The set of {\it chambers} (connected components)
of $\RR^n\backslash (\bigcup_\alpha H_\alpha)$ is
$$\cC = \{ wC \ |\ w\in W\}, \qquad\hbox{where}\qquad
C = \{ x\in \RR^n \ |\ \hbox{$\langle x,\alpha\rangle > 0$ for all $\alpha>0$}\}
$$ 
is the {\it fundamental chamber}.
A chamber $\tilde C\in \cC$ is on the {\it positive side} 
(resp. {\it negative side}) of 
the hyperplane $H_\alpha$, $\alpha\in R^+$, if 
$\langle x,\alpha\rangle>0$ (resp. $\langle x,\alpha\rangle<0$)
for all $x\in \tilde C$.
The following Proposition allows us to view
the calibration graph $\Gamma(\gamma)$ in terms of chambers in $\RR^n$.
Example 2.7 at the end of this section illustrates the conversion.

\prop Assume that $\gamma\in \RR^n$ is dominant.  Let $Z(\gamma)$ and
$P(\gamma)$ be as defined in (1.2) and let $\Gamma(\gamma)$
be the calibration graph containing $\gamma$.  

\item{(a)} The map 
$$
\matrix{
W\gamma 
&\mapleftright{1-1} 
& \left\{ \tilde C\in \cC \ |\  
\hbox{$\tilde C$ is on the positive side of
$H_\alpha$ for every $\alpha\in Z(\gamma)$} \right\} 
\cr
\cr
w\gamma &\longleftrightarrow &w^{-1}C  \cr
}
$$
is a bijection.

\item{(b)} Two vertices $u\gamma$ and $v\gamma$ in $\Gamma(\gamma)$ are 
connected by an edge if and only if the corresponding chambers $u^{-1}C$ and
$v^{-1}C$ share a face and the hyperplane
$H_{\alpha}$ containing this face satisfies
$\alpha\not\in P(\gamma)$.

\item{(c)}  The map
$$
\matrix{
\cF^{(\gamma,J)}  
&\quad \mapleftright{1-1} 
&\left\{ \tilde C\in \cC \ \ \Bigg|\ \ 
   \matrix{
   \hbox{$\tilde C$ is on the positive side of $H_\alpha$ for $\alpha\in
   Z(\gamma)$,}\hfil\cr
   \hbox{$\tilde C$ is on the positive side of $H_\alpha$ for 
   $\alpha\in P(\gamma)\backslash J$,}\hfil \cr
   \hbox{$\tilde C$ is on the negative side of $H_\alpha$ for 
               $\alpha\in J$,} \hfil\cr
   }
\right\}
\cr
\cr
w &\longleftrightarrow &w^{-1}C \cr
}
$$
is a bijection.
\pf
(a) Since $\gamma$ is dominant the stabilizer $W_\gamma$ of 
$\gamma$ is a parabolic subgroup of $W$. It is generated by
the simple reflections $s_i$ for $\alpha_i\in Z(\gamma)$.
The statement will follow from the following identifications:
$$\matrix{
W\gamma & \mapleftright{1-1} & W/W_\gamma \hfill \cr
& \mapleftright{1-1} &
\{\hbox{minimal length coset representatives of } W/W_\gamma \} \hfill \cr
& \mapleftright{1-1} &
\left\{ wC\ |\  w \hbox{ minimal length coset representative of }
W/W_\gamma\right\} \hfill \cr
& \mapleftright{1-1} &
\left\{ wC\ |\  \ell(ws_i)>\ell(w) \hbox{ for every } \alpha_i\in
Z(\gamma) \right\} \hfill \cr
& \mapleftright{1-1} &
\left\{ wC\ |\  w^{-1}C \hbox{ is on the positive side of } H_\alpha
\hbox{ for every } \alpha\in Z(\gamma) \right\} \hfill \cr
& \mapleftright{1-1} &
\left\{ w^{-1}C\ |\  w^{-1}C \hbox{ is on the positive side of } H_\alpha
\hbox{ for every } \alpha\in Z(\gamma) \right\}, \hfill \cr
}
$$
where $C = \{ x\in E\ |\ \hbox{$\langle x,\alpha\rangle >0$ for all
$\alpha>0$}\}$ is the fundamental chamber. 
Let us give a step by step explanation.
\item{(1)} There is a natural bijection between $W\gamma$ and $W/W_\gamma$.
\item{(2)} Since each coset in $W/W_\gamma$ has a unique coset representative of
minimal length (see [Bou, IV.1, Exercise~3] 
or [Hum, Proposition~1.10c]) we may identify
$W\gamma$ with the set of minimal length coset representatives of
$W/W_\gamma$. 
\item{(3)} The chambers of $\RR^n\backslash (\cup_\alpha H_\alpha)$ are 
the regions $wC, w\in W$ and thus we can identify $W\gamma$ with the set
$
\left\{ wC \ |\  
\hbox{ $w$ is a minimal length coset representative of $W/W_\gamma$}
\right\}.
$
\item{(4)} An element $w\in W$ is a minimal length coset
representative of $W/W_\gamma$ if $\ell(ws_i)>\ell(w)$ for every $i$ such
that $\alpha_i\in Z(\gamma)$. 
\item{(5)} So $w$ is a minimal length coset representative if
and only if, for any $x\in C$, $\langle w^{-1}x,\alpha_i\rangle >0$ for
every $\alpha_i\in Z(\gamma)$. 
Since $Z(\gamma)\cup -Z(\gamma)$ is a root subsystem generated by the simple roots
that it contains we have that
$w$ is a minimal length coset representative if
and only if, for any $x\in C$, $\langle w^{-1}x,\alpha\rangle >0$ for
every $\alpha\in Z(\gamma)$. 
Thus, we may identify $W\gamma$ with the set 
$
\left\{ wC \ |\ 
\hbox{$w^{-1}C$ is on the positive side of $H_\alpha$ 
for every $\alpha\in Z(\gamma)$}
\right\}.
$
\item{(6)}  The final step is to replace each chamber $wC$ by 
the chamber $w^{-1}C$.

\medskip
(b) 
Recall that $u\gamma$ and $v\gamma$ are connected by an edge in 
$\Gamma(\gamma)$ if $u\gamma=s_iv\gamma$ and 
$\langle v\gamma,\alpha_i\rangle\neq \pm 1$. If
$u$ and $v$ are minimal length coset representatives of $W/W_\gamma$,
then the chambers $u^{-1}C$ and $v^{-1}C$ should be connected by an edge if
$u=s_iv$ and $\langle \gamma,v^{-1}\alpha_i\rangle\neq\pm 1$. 
The condition that $u=s_iv$ means that $u^{-1}C=v^{-1}s_iC$.  Since 
$C$ and $s_iC$ share a face, it follows that $u^{-1}C=v^{-1}s_iC$ and
$v^{-1}C$ share a face.  This face is contained in the hyperplane
$H_{v^{-1}\alpha_i}=H_{u^{-1}\alpha_i}$ since the face that 
$C$ and $s_iC$ share is contained in the hyperplane $H_{\alpha_i}$.  
Finally, the condition
$\langle\gamma,v^{-1}\alpha_i\rangle\neq \pm 1$ is the same as
saying that $v^{-1}\alpha_i\not\in P(\gamma)$ and $u^{-1}\alpha_i\not\in
P(\gamma)$.

(c)
If $w\in W$, $\alpha>0$, $x\in C$ and $w\alpha<0$ then 
$\langle w^{-1}x,\alpha\rangle = \langle x,w\alpha\rangle <0$.
Thus
$$
R(w) = \left\{\alpha>0 \ |\  
\hbox{$w^{-1}C$ is on the negative side of $H_\alpha$ } \right\}.
$$
So the condition $R(w)\cap Z(\gamma)=\emptyset$ is equivalent 
to the condition that $w^{-1}C$ is on the positive side of $H_\alpha$ for
all $\alpha\in Z(\gamma)$. Similarly the condition $R(w)\cap P(\gamma)=J$ is
equivalent to the condition that $w^{-1}C$ be on the negative side of $H_\alpha$
for all $\alpha\in J$ and on the positive side of $H_\alpha$ for all $\alpha\in
P(\gamma)\setminus J$.
\endpf

Parts (a) and (b) of Proposition 2.2 allow us to view the calibration
graph in terms of chambers in 
$\RR^n\setminus\left(\bigcup_\alpha H_\alpha\right)$.
The vertices correspond to the chambers on the positive side
of the hyperplanes $H_\alpha$, $\alpha\in Z(\gamma)$.  The edges
of the graph are now the walls between the chambers.  The only
time that a wall between two chambers does not form an ``edge of the graph''
connecting the two chambers (vertices) is when that wall is contained
in a hyperplane $H_\alpha$ with $\alpha\in P(\gamma)$.

\subsection{Proof of Theorem 2.1}

From parts (a) and (b) of Proposition 2.2 we get that the
connected components of the graph $\Gamma(\gamma)$ correspond to the
connected components of the intersection of 
$$\left\{ \tilde C \ |\ 
\hbox{$\tilde C$ is on the positive side of $H_\alpha$ for 
$\alpha\in Z(\gamma)$} \right\} 
\qquad\hbox{with}\qquad
\RR^n\setminus \left(\bigcup_{\alpha\in P(\gamma)} H_\alpha\right).$$  
Proposition 2.2 (c) says that the standard tableaux of shape $(\gamma,J)$
correspond to the chambers which are on the positive side of
$H_\alpha$ for $\alpha\in Z(\gamma)\cup(P(\gamma)\setminus J)$
and on the negative side of $H_\alpha$ for $\alpha\in J$.  
The points in the chambers which satisfy these conditions form a 
connected subset of 
$\RR^n\setminus \left(\bigcup_{\alpha\in P(\gamma)} H_\alpha\right)$ since the
conditions describe them as the points in the intersection of
half-spaces in $\RR^n$.
Finally, one only needs to note that the connected components
$\RR^n\setminus \left(\bigcup_{\alpha\in P(\gamma)} H_\alpha\right)$  
are determined by which subset of hyperplanes in $P(\gamma)$ 
they are on the negative side of.
This completes the proof of Theorem 2.1.
\endpf

\subsection{Invariance properties for calibration graphs.}

\smallskip\noindent
{\sl Invariance property (1):}
As unlabeled graphs, $\Gamma(\gamma)=\Gamma(w\gamma)$ for
all $w\in W$. 
\smallskip\noindent
{\sl Invariance property (2):}  
If $Z(\gamma)=Z(\kappa)$ and $P(\gamma)=P(\kappa)$
then $\Gamma(\gamma)=\Gamma(\kappa)$.
\medskip\noindent
Property (1) is immediate from the definition of $\Gamma(\gamma)$.
From the definition of the calibration graph one sees that the edges of 
$\Gamma(\gamma)$ are controlled by the set $P(\gamma)$.  
Since the vertices $W\gamma$ can be identified with the set $W/W_\gamma$, 
where $W_\gamma$ is the stabilizer of $\gamma$, it follows that
the graph $\Gamma(\gamma)$
depends only on $W_\gamma$ and the set $P(\gamma)$.  
Since $W_\gamma$ is the group generated by the reflections $s_\alpha$ for
$\alpha\in Z(\gamma)$, it
follows that the structure of $\Gamma(\gamma)$ is dependent only on the sets
$Z(\gamma)$ and $P(\gamma)$. This establishes invariance property (2).

\subsection{Intersections and shapes.}

Let $\cA$ be the arrangement of (affine) hyperplanes
given by
$$\cA' = \{ H_{\alpha-\delta}, H_\alpha, H_{\alpha+\delta} \ |\ \alpha\in R^+\},
\quad\hbox{where}\quad
\eqalign{
H_\alpha &= \{ x\in \RR^n \ |\ \langle x,\alpha\rangle=0 \}, \cr
H_{\alpha+\delta} &= \{ x\in \RR^n \ |\ \langle x,\alpha\rangle=1 \}, \cr
H_{\alpha-\delta} &= \{ x\in \RR^n \ |\ \langle x,\alpha\rangle=-1 \}. \cr
}$$
The intersection lattice $L(\cA')$ is the set of {\it intersections},
$I=\bigcap_{H_\beta\in \cB} H_\beta$,
$\cB\subseteq \cA',$
partially ordered by inclusion (as subsets of $\RR^n$).
Since $\cA'$ is symmetric under the Weyl group $L(\cA')$
also carries a Weyl group symmetry.  The quotient $L(\cA')/W$ is
constructed by  identifying intersections $I_1$ and $I_2$ if there is 
a $w\in W$ such that $wI_1=I_2$.  
It follows from the invariance properties of the calibration graphs
that the distinct calibration graphs
are in one-to-one correspondence with the elements of the quotient $L(\cA')/W$.
In particular, the number of distinct calibration graphs is {\it finite}.

Let $(\gamma,J)$ and $(\eta,K)$ be two placed shapes.  We shall say that 
$(\gamma,J)$ and $(\eta,K)$ have the same underlying {\it shape} if
there is a $w\in W$ such that
$$Z(w\gamma)=Z(\eta), \qquad P(w\gamma)=P(\eta),
\qquad\hbox{and}\qquad
wJ=K.$$
When $(\gamma,J)$ and $(\eta,K)$ are placed shapes with the same
underlying shape then the calibration graphs $\Gamma(\gamma)$ and 
$\Gamma(\kappa)$ are the same and there is 
a natural bijection
$$\cF^{(\gamma,J)} ~\mapleftright{1-1}~ \cF^{(\eta,K)}.$$

\subsection{}

It would be interesting to determine the size of
$L(\cA')/W$, i.e. the number of distinct calibration graphs.
A. Postnikov has explained that in the Type A case the numbers
$f(n)={\rm Card}(L(\cA')/S_n)$ have the following generating function
$$\prod_{k\ge 1} (1-q^k)^{2^{k-1}} = \sum_{n\ge 1} f(n)q^n.$$
One can prove this by counting the number of ways of constructing
the configurations of boxes described in Section 3.

The arrangement $\cA'$ is very similar to the Shi arrangement
$$\cA'' = \{ H_\alpha, H_{\alpha+\delta} \ |\ \alpha\in R^+\}.$$
C. Athanasiadis has told me that, for the type $A$ case,
the number of elements of the intersection lattice $L(\cA'')$ which contain a
dominant weight is equal to
$$\sum_{k=1}^n {n-1 \choose k-1} F_{2k-1}$$
where $F_1=1, \ F_3=2, \ F_5=5, \ldots$ are the odd Fibonacci numbers.
The Shi arrangement has been an object of intense
recent study, see [Sh3], [St1], [AL].  
There are many indications
[Sh3], [Xi,1.11, 2.6] that there should be a strong connection between the
work in this paper and the combinatorics of the Shi arrangement.

\subsection{Example.}

Consider the root system of type $C_2$ where
$R^+ = \{ \alpha_1,\alpha_2,\alpha_1+\alpha_2,\alpha_1+2\alpha_2\}$.
Realize this root system in $\RR^2$ by letting
$\alpha_1=2\varepsilon_1$ and $\alpha_2=\varepsilon_2-\varepsilon_1$,
where $\{\varepsilon_1,\varepsilon_2\}$ is an orthonormal basis of
$\RR^2$.  Let $\gamma\in \RR^2$ be given by
$\langle \gamma,\alpha_1\rangle = 0$
and
$\langle \gamma,\alpha_2\rangle =1$.
Then $\gamma$ is dominant (i.e. in $\overline{C}$) and integral and
$$Z(\gamma)=\{\alpha_1\}
\qquad\hbox{and}\qquad
P(\gamma)=\{\alpha_2,\alpha_1+\alpha_2\}.$$
$$
\matrix{
\beginpicture
\setcoordinatesystem units <1cm,1cm>         % sets scale
\setplotarea x from -4 to 4, y from -4 to 4    % sets plot size up
\put{$H_{\alpha_1}$}[b] at 0 3.1
\put{$H_{\alpha_2}$}[l] at 3.1 3.1
\put{$H_{\alpha_1+\alpha_2}$}[r] at -3.1 3.1
\put{$H_{\alpha_1+2\alpha_2}$}[l] at 3.1 0
\put{$H_{\alpha_1+2\alpha_2+\delta}$}[tr] at -3.1 0.9
\put{$H_{\alpha_1+2\alpha_2-\delta}$}[br] at -3.1 -0.9
\put{$H_{\alpha_1+\alpha_2+\delta}$}[tl] at 3.1 -1.6
\put{$H_{\alpha_1+\alpha_2-\delta}$}[tl] at 1.6 -3.5
\put{$H_{\alpha_2+\delta}$}[tr] at -3.5 -1.7
\put{$H_{\alpha_2-\delta}$}[tr] at -1.6 -3.5
\put{$H_{\alpha_1-\delta}$}[t] at -0.9 -3.5
\put{$H_{\alpha_1+\delta}$}[t] at  0.9 -3.5
\put{$\bullet$} at 0 2
\put{$\gamma$}[r] at -0.1 2.1
\put{$\bullet$} at 2 0
\put{$s_2\gamma$}[tl] at 2.05 -0.1
\put{$\bullet$} at -2 0 
\put{$s_1s_2\gamma$}[tl] at -3 -0.1
\put{$\bullet$} at 0 -2
\put{$s_2s_1s_2\gamma$}[tl] at .04 -2.04 
\plot -3 -3   3 3 /
\plot  3 -3  -3 3 /
\plot  0  3   0 -3 /
\plot  3  0  -3  0 /
\setdots
\putrule from 1 3.5 to 1 -3.5
\putrule from -1 3.5 to -1 -3.5
\putrule from 3.5 1 to -3.5  1
\putrule from 3.5 -1 to -3.5 -1
\plot 3.5 1.5  -1.5 -3.5 /
\plot -3.5 -1.5  1.5 3.5 /
\plot -3.5 1.5  1.5 -3.5 /
\plot 3.5 -1.5  -1.5 3.5 /
\endpicture
\cr
\cr
\hbox{The arrangement $\cA'$}
\cr
}$$
$$\matrix{
\beginpicture
\setcoordinatesystem units <1cm,1cm>         % sets scale
\setplotarea x from -3 to 5, y from -3 to 3    % sets plot size up
\put{$J=\emptyset$}[l] at 4 2
\put{$J=\{\alpha_2\}$}[l] at 4 0
\put{$J=\{\alpha_2,\alpha_1+\alpha_2\}$}[l] at 4 -2 
\put{$\bullet$} at 0 2
\put{$\gamma$}[r] at -0.1 2.1
\put{$\bullet$} at 2 0
\put{$s_2\gamma$}[tl] at 2.1 -0.1
\put{$\bullet$} at -2 0 
\put{$s_1s_2\gamma$}[tr] at -2.1 -0.1
\put{$\bullet$} at 0 -2
\put{$s_2s_1s_2\gamma$}[tl] at .01 -2.1 
\plot  -2 0   2 0 /
\endpicture
\cr
\cr
\hbox{The calibration graph $\Gamma(\gamma)$} \cr
}
$$
When we convert to regions in 
$\RR^2\setminus \left(\bigcup_\alpha H_\alpha\right)$
we get the picture
$$
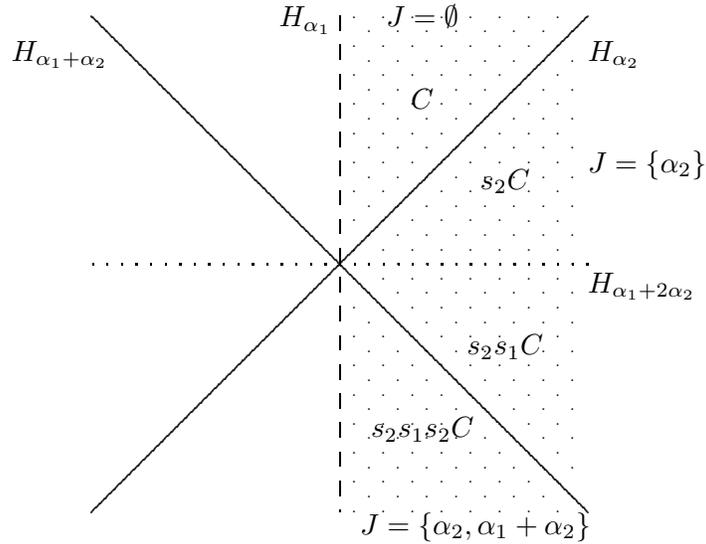
\beginpicture
\setcoordinatesystem units <1cm,1cm>         % sets scale
\setplotarea x from -3 to 3, y from -3 to 3    % sets plot size up
\put{$H_{\alpha_1}$}[r] at -0.1 2.95
\put{$H_{\alpha_2}$}[tl] at 3 2.7
\put{$H_{\alpha_1+\alpha_2}$}[tr] at -2.8 2.7
\put{$H_{\alpha_1+2\alpha_2}$}[tl] at 3 -0.1
\put{$J=\emptyset$} at 1 3
\put{$J=\{\alpha_2\}$}[l] at 3 1.2
\put{$J=\{\alpha_2,\alpha_1+\alpha_2\}$}[tl] at 0.25 -3.0
\put{$C$} at 1 2
\put{$s_2C$} at 2 1
\put{$s_2s_1C$} at 2 -1
\put{$s_2s_1s_2C$} at  1 -2
\plot -3 -3   3 3 /
\plot  3 -3  -3 3 /
\setdashes
\putrule from 0 3 to 0 -3
\setdots
\putrule from 3 0 to -3  0
\vshade 0 -3 3    3 -3 3 /
\endpicture
$$
The dashed line is the hyperplane corresponding to the root in 
$Z(\gamma)$ and the solid lines are the hyperplanes 
corresponding to the roots in $P(\gamma)$.
\endexample

\section 3. Type {\rm A} and configurations of boxes

\subsection{The root system.}

Let $\{\varepsilon_1,\ldots,\varepsilon_n\}$ be an
orthonormal basis  of
$\RR^n$ so that each sequence $\gamma=(\gamma_1,\ldots,\gamma_n)\in \RR^n$
is identified with the vector $\gamma=\sum_i \gamma_i\varepsilon_i$.
The root system of type $A_{n-1}$ is given by the sets
$$
R = \left\{\varepsilon_j-\varepsilon_i\,|\, 1\leq i,j\leq n\right\}
\quad\hbox{and}\quad
R^+=\left\{\varepsilon_j-\varepsilon_i\,|\, 1\leq i<j\leq n\right\}.
$$
The Weyl group is $W=S_n$, the symmetric group, acting by 
permutations of the $\varepsilon_i$.

\subsection{Partitions, skew shapes, and (classical) standard tableaux.}

A partition $\lambda$ is a collection of $n$
boxes in a  corner.  We shall conform to the conventions in [Mac] and assume
that gravity goes up and to the left. 

$$
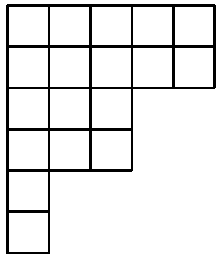
\beginpicture
\setcoordinatesystem units <0.5cm,0.5cm>         % sets scale
\setplotarea x from 0 to 4, y from 0 to 3    % sets plot size up
\linethickness=0.5pt                          % sets line thickness
\putrule from 0 6 to 5 6          %            
\putrule from 0 5 to 5 5          %  draws horizontal lines         
\putrule from 0 4 to 5 4          %           
\putrule from 0 3 to 3 3          %           
\putrule from 0 2 to 3 2          %           
\putrule from 0 1 to 1 1          %           
\putrule from 0 0 to 1 0          %  

\putrule from 0 0 to 0 6        %           
\putrule from 1 0 to 1 6        %           
\putrule from 2 2 to 2 6        %             
\putrule from 3 2 to 3 6        %  draws vertical lines           
\putrule from 4 4 to 4 6        %             
\putrule from 5 4 to 5 6        %  
%\vshade 2 1 2   3 1 2 /           
\endpicture
$$
Any partition $\lambda$ can be identified with
the sequence $\lambda=(\lambda_1\ge \lambda_2\ge \ldots )$
where $\lambda_i$ is the number of boxes in row $i$ of $\lambda$.
The rows and columns are numbered in the same way as for matrices.
In the example above we have $\lambda=(553311)$.

If $\lambda$ and $\mu$ are partitions such that $\mu_i\le \lambda_i$
for all $i$ we write $\mu\subseteq \lambda$.  
The {\it skew shape} $\lambda/\mu$
consists of all boxes of $\lambda$ which are not in $\mu$.
Let $\lambda/\mu$ be a skew shape with $n$ boxes.
Number the boxes of each skew shape
$\lambda/\mu$ along major diagonals from southwest to northeast and
$$\hbox{write ${\rm box}_i$ to indicate the box numbered $i$.}$$
See Example 3.4  below.
A (classical) {\it standard tableau of shape} $\lambda/\mu$ is a filling  
of the boxes in the skew shape
$\lambda/\mu$ with the numbers $1,\ldots,n$ such that 
the numbers increase from left to right in each row and from top to
bottom down each column. 
Let $\cF^{\lambda/\mu}$ be the set of standard tableaux of shape
$\lambda/\mu$.
Given a standard tableau $t$ of shape $\lambda/\mu$ define the 
{\it word} of $t$ to be permutation
$$
w_t=\left(\matrix{
1&\cdots & n \cr
t({\rm box}_1)&\ldots &t({\rm box}_n)
}\right)
$$ 
where $t({\rm box}_i)$ is the entry in ${\rm box}_i$ of $t$. 

\subsection{Placed skew shapes}

Let $\lambda/\mu$ be a skew shape with $n$ boxes.  Imagine placing
$\lambda/\mu$ on a piece of infinite graph paper where the diagonals of the 
graph paper are indexed consecutively (with elements of $\ZZ$) from southeast 
to northwest.  The {\it content} of a box $b$ is the number $c(b)$ of the
diagonal that $b$ is on.
Let 
$$\gamma=\sum_{i=1}^n c({\rm box}_i)\varepsilon_i,$$
which we identify with the sequence 
$\gamma=(c({\rm box}_1),c({\rm box}_2),\ldots,c({\rm box}_n))$.
The pair $(\gamma,\lambda/\mu)$ is a {\it placed skew shape}.
It follows from the definitions in (1.2) that
$$
\eqalign{
Z(\gamma) &= \{ \varepsilon_j-\varepsilon_i \ |\ 
\hbox{$j>i$ and 
${\rm box}_j$ and ${\rm box}_i$ are in the same diagonal} \}, 
\quad\hbox{and} \cr
P(\gamma) &= \{ \varepsilon_j-\varepsilon_i \ |\ 
\hbox{$j>i$ and 
${\rm box}_j$ and ${\rm box}_i$ are in adjacent diagonals} \}. \cr
}
$$
Define 
$$
J=\left\{\varepsilon_j-\varepsilon_i\ \  \Bigg|\ \ 
\matrix{
j>i \hfill \cr 
\hbox{${\rm box}_j$ and ${\rm box}_i$ are in adjacent diagonals} \hfill \cr
\hbox{${\rm box}_j$ is northwest of ${\rm box}_i$ } \hfill\cr
} 
\right\},
$$
where ${\it northwest}$ means strictly north and weakly west.

\subsection{Example.}

The following diagrams illustrate standard tableaux and
the numbering of boxes in a skew shape $\lambda/\mu$.
$$
\matrix{
\beginpicture
\setcoordinatesystem units <0.5cm,0.5cm>         % sets scale
\setplotarea x from 0 to 4, y from 0 to 3    % sets plot size up
\linethickness=0.5pt                          % sets line thickness
\put{1} at 0.5 0.5
\put{2} at 0.5 1.5
\put{3} at 1.5 1.5
\put{4}  at 3.5 2.5
\put{5} at 4.5 3.5
\put{6} at 4.5 4.5
\put{7}  at 5.5 3.5
\put{8}  at 5.5 4.5
\put{10}  at 5.5 5.5
\put{9}  at 6.5 3.5
\put{11}  at 6.5 4.5
\put{12}  at 6.5 5.5
\put{13}  at 7.5 5.5
\put{14}  at 8.5 5.5
\putrule from 5 6 to 9 6          %            
\putrule from 4 5 to 9 5          %  draws horizontal lines         
\putrule from 4 4 to 7 4          %           
\putrule from 3 3 to 7 3          %           
\putrule from 3 2 to 4 2          %           
\putrule from 0 2 to 2 2          %           
\putrule from 0 1 to 2 1          %           
\putrule from 0 0 to 1 0          %  
\putrule from 0 0 to 0 2        %           
\putrule from 1 0 to 1 2        %           
\putrule from 2 1 to 2 2        %             
\putrule from 3 2 to 3 3        %  draws vertical lines           
\putrule from 4 2 to 4 5        %             
\putrule from 5 3 to 5 6        %  
\putrule from 6 3 to 6 6        %  
\putrule from 7 3 to 7 6        %  
\putrule from 8 5 to 8 6        %  
\putrule from 9 5 to 9 6        %  
%\vshade 2 1 2   3 1 2 /           
\endpicture
&\qquad\qquad
&
\beginpicture
\setcoordinatesystem units <0.5cm,0.5cm>         % sets scale
\setplotarea x from 0 to 4, y from 0 to 3    % sets plot size up
\linethickness=0.5pt                          % sets line thickness
\put{11} at 0.5 0.5
\put{6} at 0.5 1.5
\put{8} at 1.5 1.5
\put{2}  at 3.5 2.5
\put{7} at 4.5 3.5
\put{1} at 4.5 4.5
\put{13}  at 5.5 3.5
\put{5}  at 5.5 4.5
\put{3}  at 5.5 5.5
\put{14}  at 6.5 3.5
\put{10}  at 6.5 4.5
\put{4}  at 6.5 5.5
\put{9}  at 7.5 5.5
\put{12}  at 8.5 5.5
\putrule from 5 6 to 9 6          %            
\putrule from 4 5 to 9 5          %  draws horizontal lines         
\putrule from 4 4 to 7 4          %           
\putrule from 3 3 to 7 3          %           
\putrule from 3 2 to 4 2          %           
\putrule from 0 2 to 2 2          %           
\putrule from 0 1 to 2 1          %           
\putrule from 0 0 to 1 0          %  
\putrule from 0 0 to 0 2        %           
\putrule from 1 0 to 1 2        %           
\putrule from 2 1 to 2 2        %             
\putrule from 3 2 to 3 3        %  draws vertical lines           
\putrule from 4 2 to 4 5        %             
\putrule from 5 3 to 5 6        %  
\putrule from 6 3 to 6 6        %  
\putrule from 7 3 to 7 6        %  
\putrule from 8 5 to 8 6        %  
\putrule from 9 5 to 9 6        %  
%\vshade 2 1 2   3 1 2 /           
\endpicture
\cr
\hbox{$\lambda/\mu$ with boxes numbered}
&&\hbox{A standard tableau $t$ of shape $\lambda/\mu$} \cr
}
$$
The word of the standard tableau $t$ is the permutation
$w_t=(11,6,8,2,7,1,13,5,14,3,10,4,9,12)$ (in one-line 
notation).

The following picture shows the contents of the boxes in the
placed skew shape $(\gamma,\lambda/\mu)$ with
$\gamma = (-7,-6,-5,-2,0,1,1,2,2,3,3,4,5,6)$.
$$
\matrix{
\beginpicture
\setcoordinatesystem units <0.5cm,0.5cm>         % sets scale
\setplotarea x from 0 to 4, y from 0 to 3    % sets plot size up
\linethickness=0.5pt                          % sets line thickness
\put{-7} at 0.5 0.5
\put{-6} at 0.5 1.5
\put{-5} at 1.5 1.5
\put{-2}  at 3.5 2.5
\put{0} at 4.5 3.5
\put{1} at 4.5 4.5
\put{1}  at 5.5 3.5
\put{2}  at 5.5 4.5
\put{3}  at 5.5 5.5
\put{2}  at 6.5 3.5
\put{3}  at 6.5 4.5
\put{4}  at 6.5 5.5
\put{5}  at 7.5 5.5
\put{6}  at 8.5 5.5
\putrule from 5 6 to 9 6          %            
\putrule from 4 5 to 9 5          %  draws horizontal lines         
\putrule from 4 4 to 7 4          %           
\putrule from 3 3 to 7 3          %           
\putrule from 3 2 to 4 2          %           
\putrule from 0 2 to 2 2          %           
\putrule from 0 1 to 2 1          %           
\putrule from 0 0 to 1 0          %  
\putrule from 0 0 to 0 2        %           
\putrule from 1 0 to 1 2        %           
\putrule from 2 1 to 2 2        %             
\putrule from 3 2 to 3 3        %  draws vertical lines           
\putrule from 4 2 to 4 5        %             
\putrule from 5 3 to 5 6        %  
\putrule from 6 3 to 6 6        %  
\putrule from 7 3 to 7 6        %  
\putrule from 8 5 to 8 6        %  
\putrule from 9 5 to 9 6        %  
%\vshade 2 1 2   3 1 2 /           
\endpicture
\cr
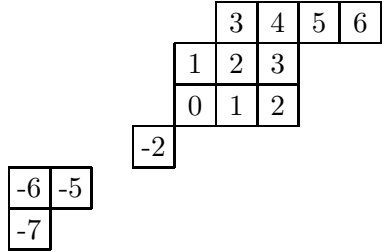
\hbox{Contents of the boxes of $(\gamma,\lambda/\mu)$} \cr
}
$$
In this case
$J=\{
\varepsilon_2-\varepsilon_1,
\varepsilon_6-\varepsilon_5,
\varepsilon_8-\varepsilon_7,
\varepsilon_{10}-\varepsilon_8,
\varepsilon_{10}-\varepsilon_9,
\varepsilon_{11}-\varepsilon_9,
\varepsilon_{12}-\varepsilon_{11}\}.$
\endpf

\subsection{}

The following theorem shows how the generalized standard tableaux
defined in (1.2) reduce to the classical standard Young tableaux
in the type $A$ case.

\thm  Let $(\gamma,\lambda/\mu)$ be a placed skew shape and
let $J$ be as defined in (3.3).  Let
$\cF^{\lambda/\mu}$ be the set of standard tableaux of shape $\lambda/\mu$
and let $\cF^{(\gamma,J)}$ be the set of generalized standard tableaux of shape
$(\gamma,J)$ as defined in (1.2). 
Then the map 
$$
\matrix{
\cF^{\lambda/\mu} &\mapleftright{1-1} & \cF^{(\gamma,J)} \cr
t & \longleftrightarrow & w_t, \cr
}
$$
where $w_t$ is as defined in (3.2),
is a bijection.
\pf
If $w=(w(1)\cdots w(n))$ is a permutation in $S_n$ then 
$$R(w)=\{\varepsilon_j-\varepsilon_i\ |\ 
\hbox{$j>i$ such that $w(j)<w(i)$} \ \}.$$  
The theorem is a consequence of the following chain of equivalences:

\medskip\noindent
The filling $t$ is a standard tableau if and only if for all $1\le i<j\le n$
\item{(a)} $t({\rm box}_i)<t({\rm box}_j)$ 
if ${\rm box}_i$ and ${\rm box}_j$ are on the same diagonal,
\item{(b)} $t({\rm box}_i)<t({\rm box}_j)$ 
if ${\rm box}_j$ is immediately to the right of ${\rm box}_i$, and
\item{(c)} $t({\rm box}_i)>t({\rm box}_j)$ 
if ${\rm box}_j$ is  immediately above ${\rm box}_i$.
\smallskip\noindent
These conditions hold if and only if
\item{(a)} $\varepsilon_j-\varepsilon_i\not\in R(w(t))$ if
$\varepsilon_j-\varepsilon_i\in Z(\gamma)$,
\item{(b)} $\varepsilon_j-\varepsilon_i\not\in R(w(t))$ if
$\varepsilon_j-\varepsilon_i\in P(\gamma)\setminus J$,
\item{(c)} $\varepsilon_j-\varepsilon_i\in R(w(t))$ if
$\varepsilon_j-\varepsilon_i\in J$,
\smallskip\noindent
which hold if and only if
$$\hbox{(a)\enspace $\alpha\not\in R(w(t))$ if $\alpha\in Z(\gamma)$,}
\quad
\hbox{(b)\enspace $\alpha\not\in R(w(t))$ if $\alpha\in
P(\gamma)\setminus J$,}
\quad\hbox{and}\quad
\hbox{(c)\enspace $\alpha\in R(w(t))$ if $\alpha\in J$.}
$$
Finally, these are equivalent to the conditions
$R(w(t))\cap Z(\gamma)=\emptyset$ and 
$R(w(t))\cap P(\gamma)=J$.
\endpf

\subsection{Placed configurations}

We have described how one can identify placed
skew shapes $(\gamma,\lambda/\mu)$ with pairs $(\gamma, J)$.
One can extend this conversion to associate placed configurations of boxes
to more general pairs $(\gamma,J)$.  The resulting configurations
are not always skew shapes.

Let $(\gamma,J)$ be a pair such that $\gamma=(\gamma_1,\ldots,\gamma_n)$
is a dominant integral weight and $J\subseteq P(\gamma)$. 
(The sequence $\gamma$ is a dominant integral weight
if $\gamma_1\leq\cdots\leq\gamma_n$ and $\gamma_i\in\ZZ$ for all $i$.)
If $J$ satisfies the condition
$$\hbox{{\sl If $\beta\in J$, $\alpha\in Z(\gamma)$, and $\beta-\alpha\in R^+$
then $\beta-\alpha\in J$}}
$$
then $(\gamma,J)$ will determine a placed configuration of boxes.
As in the placed skew shape case, think of the boxes as being placed on graph
paper where the  boxes on a given diagonal all have the same content.  
(The boxes on each diagonal are allowed to slide along the diagonal
as long as they don't pass through the corner of a box on an adjacent diagonal.)
The sequence $\gamma$ describes how many boxes are on each diagonal and the set
$J$ determines how the boxes on adjacent diagonals are placed relative to each
other.  We want
$$\gamma=\sum_{i=1}^n c({\rm box}_i)\varepsilon_i,$$
and
\itemitem{(a)} If $\varepsilon_j-\varepsilon_i\in J$ then ${\rm box}_j$ is
northwest of ${\rm box}_i$, and
\itemitem{(b)} If $\varepsilon_j-\varepsilon_i\in P(\gamma)\backslash J$ then
${\rm box}_j$ is southeast of ${\rm box}_i$,
\smallskip\noindent
where the boxes are numbered along diagonals in the same way as 
for skew shapes,
{\it southeast} means weakly south and strictly east, 
and {\it northwest} means strictly north and weakly west. 

If we view the pair $(\gamma,J)$ as a placed configuration of boxes then
the {\it standard tableaux} are fillings $t$ of the $n$ boxes in the
configuration with $1,2,\ldots, n$ such that for all $i<j$
\item{(a)} $t({\rm box}_i)<t({\rm box}_j)$ if ${\rm box}_i$ and ${\rm
box}_j$ are on the same diagonal,
\item{(b)} $t({\rm box}_i)<t({\rm box}_j)$ if ${\rm box}_i$ and ${\rm
box}_j$ are on adjacent diagonals and ${\rm box}_j$ is southeast 
of ${\rm box}_i$, and
\item{(c)} $t({\rm box}_i)>t({\rm box}_j)$ if ${\rm box}_i$ and ${\rm 
box}_j$ are on adjacent diagonals and ${\rm box}_j$ is northwest 
of ${\rm box}_i$.
\smallskip\noindent
As in Theorem 3.5 the permutation in $\cF^{(\gamma,J)}$ which corresponds
to the standard tableau $t$ is
$w(t)=(t({\rm box}_1),\ldots,t({\rm box}_n))$.
The following example illustrates the conversion.

\bigskip\noindent
{\sl Example.}
Suppose $\gamma=(-1,-1,-1,0,0,0,1,1,1,2,2,2)$ and
$$\eqalign{
J &=\left\{ \varepsilon_4-\varepsilon_1, \varepsilon_4-\varepsilon_2,
\varepsilon_4-\varepsilon_3, \varepsilon_5-\varepsilon_2,
\varepsilon_5-\varepsilon_3, \varepsilon_7-\varepsilon_5,
\varepsilon_7-\varepsilon_6, \varepsilon_8-\varepsilon_6,
\varepsilon_{10}-\varepsilon_9, \varepsilon_{10}-\varepsilon_8,\right. \cr
&\phantom{J = } \left.\varepsilon_{10}-\varepsilon_7,
\varepsilon_{11}-\varepsilon_9,
\varepsilon_{11}-\varepsilon_8, \varepsilon_{11}-\varepsilon_7,
\varepsilon_{12}-\varepsilon_9 \right\}. \cr
}
$$
The placed configuration of boxes corresponding to $(\gamma,J)$
is as given below.
$$
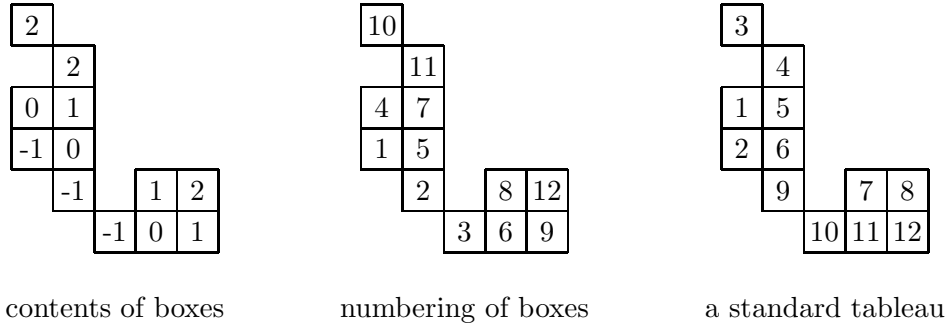
\matrix{
\beginpicture
\setcoordinatesystem units <0.5cm,0.5cm> point at 6 0       % sets scale
\setplotarea x from 0 to 5, y from -1 to 6    % sets plot size up
\linethickness=0.5pt                          % sets line thickness
%\put{contents of boxes}  at 2.5 -0.5    % 
\put{-1}  at 0.5 2.5    % 
\put{-1}  at 1.5 1.5    % 
\put{-1}  at 2.5 0.5    % 
\put{0}  at 0.5 3.5    % 
\put{0}  at 1.5 2.5    %  labels contents of boxes
\put{0}  at 3.5 0.5    % 
\put{1}  at 1.5 3.5    % 
\put{1}  at 3.5 1.5    % 
\put{1}  at 4.5 0.5    % 
\put{2}  at 0.5 5.5    % 
\put{2}  at 1.5 4.5    % 
\put{2}  at 4.5 1.5    % 
\putrule from 0 6 to 1 6          %            
\putrule from 0 5 to 2 5          %            
\putrule from 0 4 to 2 4          %            
\putrule from 0 3 to 2 3          %            
\putrule from 3 2 to 5 2          %        
\putrule from 0 2 to 2 2          %  draws horizontal lines         
\putrule from 1 1 to 5 1          %           
\putrule from 2 0 to 5 0          %           
\putrule from 0 2 to 0 4        %             
\putrule from 0 5 to 0 6        %             
\putrule from 1 1 to 1 6        %             
\putrule from 2 0 to 2 5        %  draws vertical lines           
\putrule from 3 0 to 3 2        %             
\putrule from 4 0 to 4 2        %  
\putrule from 5 0 to 5 2        %  
\endpicture
&\qquad
&\beginpicture
\setcoordinatesystem units <0.5cm,0.5cm> point at -1 0       % sets scale
\setplotarea x from 0 to 5, y from -1 to 6    % sets plot size up
\linethickness=0.5pt                          % sets line thickness
%\put{numbering of boxes}  at 2.5 -0.5    % 
\put{1}  at 0.5 2.5    % 
\put{2}  at 1.5 1.5    % 
\put{3}  at 2.5 0.5    % 
\put{4}  at 0.5 3.5    % 
\put{5}  at 1.5 2.5    %  labels numbers of boxes
\put{6}  at 3.5 0.5    % 
\put{7}  at 1.5 3.5    % 
\put{8}  at 3.5 1.5    % 
\put{9}  at 4.5 0.5    % 
\put{10} at 0.5 5.5    % 
\put{11} at 1.5 4.5    % 
\put{12} at 4.5 1.5    % 
\putrule from 0 6 to 1 6          %            
\putrule from 0 5 to 2 5          %            
\putrule from 0 4 to 2 4          %            
\putrule from 0 3 to 2 3          %            
\putrule from 3 2 to 5 2          %        
\putrule from 0 2 to 2 2          %  draws horizontal lines         
\putrule from 1 1 to 5 1          %           
\putrule from 2 0 to 5 0          %           
\putrule from 0 2 to 0 4        %             
\putrule from 0 5 to 0 6        %             
\putrule from 1 1 to 1 6        %             
\putrule from 2 0 to 2 5        %  draws vertical lines           
\putrule from 3 0 to 3 2        %             
\putrule from 4 0 to 4 2        %  
\putrule from 5 0 to 5 2        %  
\endpicture 
&\qquad
&\beginpicture
\setcoordinatesystem units <0.5cm,0.5cm> point at -1 0       % sets scale
\setplotarea x from 0 to 5, y from -1 to 6    % sets plot size up
\linethickness=0.5pt                          % sets line thickness
%\put{numbering of boxes}  at 2.5 -0.5    % 
\put{2}  at 0.5 2.5    % 
\put{9}  at 1.5 1.5    % 
\put{10}  at 2.5 0.5    % 
\put{1}  at 0.5 3.5    % 
\put{6}  at 1.5 2.5    %  labels numbers of boxes
\put{11}  at 3.5 0.5    % 
\put{5}  at 1.5 3.5    % 
\put{7}  at 3.5 1.5    % 
\put{12}  at 4.5 0.5    % 
\put{3} at 0.5 5.5    % 
\put{4} at 1.5 4.5    % 
\put{8} at 4.5 1.5    % 
\putrule from 0 6 to 1 6          %            
\putrule from 0 5 to 2 5          %            
\putrule from 0 4 to 2 4          %            
\putrule from 0 3 to 2 3          %            
\putrule from 3 2 to 5 2          %        
\putrule from 0 2 to 2 2          %  draws horizontal lines         
\putrule from 1 1 to 5 1          %           
\putrule from 2 0 to 5 0          %           
\putrule from 0 2 to 0 4        %             
\putrule from 0 5 to 0 6        %             
\putrule from 1 1 to 1 6        %             
\putrule from 2 0 to 2 5        %  draws vertical lines           
\putrule from 3 0 to 3 2        %             
\putrule from 4 0 to 4 2        %  
\putrule from 5 0 to 5 2        %  
\endpicture 
\cr
\hbox{contents of boxes}
&&\hbox{numbering of boxes}
&&\hbox{a standard tableau} \cr
}
$$
\endexample

\subsection{Books of placed configurations.}

The most general case to consider is when
$\gamma=(\gamma_1,\ldots,\gamma_n)$ is an arbitrary element of
$\RR^n$ and $J\subseteq P(\gamma)$.  This case is handled as follows.
First group the entries of $\gamma$ according to their $\ZZ$-coset in $\RR$.
Each group of entries in $\gamma$ can be arranged to form a sequence
$$
\beta+C_\beta = \beta+(z_1,\ldots,z_k) = (\beta+z_1,\ldots,\beta+z_k),
\quad\hbox{
where $0\leq\beta<1$, $z_i\in\ZZ$ and $z_1\leq \cdots
\leq z_k$.}
$$
Fix some ordering of these groups and let 
$$\vec{\gamma}=(\beta_1+C_{\beta_1},\ldots,\beta_r+C_{\beta_r})$$
be the rearrangement of the sequence $\gamma$ with the groups listed in
order.  Since $\vec{\gamma}$ is a rearrangement of $\gamma$,
the calibration graphs $\Gamma(\gamma)$ and $\Gamma(\vec{\gamma})$
are the same (see Invariance Property (1) in (2.4)).  This means that
it is sufficient to understand the standard tableaux corresponding
to $\vec\gamma$. 

The decomposition of $\vec \gamma$ into
groups induces decompositions
$$
Z(\vec{\gamma})=\bigcup_{\beta_i}Z_{\beta_i}, \quad\quad
P(\vec{\gamma})=\bigcup_{\beta_i}P_{\beta_i}, \quad\hbox{
and, if $J\subseteq P(\vec{\gamma})$, then }\quad
J=\bigcup_{\beta_i}J_{\beta_i},$$
where $J_{\beta_i}=J\cap P_{\beta_i}$.
Each pair $(C_\beta,J_{\beta})$ is a placed shape of the type considered
in the previous subsection and we may
identify $(\vec{\gamma},J)$ with the {\it book of placed shapes}
$\left( (C_{\beta_1},J_{\beta_1}),\ldots,(C_{\beta_r},J_{\beta_r}) \right)$.
We think of this as a {\it book} with {\it pages} numbered by the values
$\beta_1,\ldots,\beta_r$ and with the placed configuration determined by
$(C_{\beta_i},J_{\beta_i})$ on page $\beta_i$.
In this form the {\it standard tableaux} of shape $(\vec \gamma,J)$ are
fillings of the $n$ boxes in the book with the numbers $1,\ldots,n$ such
that the filling on each page satisfies the conditions for a standard
tableau in (3.6).
 
\bigskip\noindent
{\sl Example.}   
If
$\gamma=(1/2,1/2,1,1,1,3/2,-2,-2,-1/2,-1,-1,-1,-1/2,1/2,0,0,0)$ then one
possibility for $\vec{\gamma}$ is
$$
\vec{\gamma} = (-2,-2,-1,-1,-1,0,0,0,1,1,1,-1/2,-1/2,1/2,1/2,1/2,3/2).$$
In this case $\beta_1=0$, $\beta_2=1/2$,
$$\beta_1+C_{\beta_1}=(-2,-2,-1,-1,-1,0,0,0,1,1,1)\quad\hbox{and}\quad
\beta_2+C_{\beta_2}=(-1/2,-1/2,1/2,1/2,1/2,3/2).$$
If $J=J_{\beta_1}\cup J_{\beta_2}$ where 
$J_{\beta_2} = \{ \varepsilon_{14}-\varepsilon_{13}, 
\varepsilon_{17}-\varepsilon_{16} \}$ and 
$$
J_{\beta_1} = \{ \varepsilon_3-\varepsilon_2, \varepsilon_4-\varepsilon_2,
\varepsilon_5-\varepsilon_2, \varepsilon_6-\varepsilon_3,
\varepsilon_6-\varepsilon_4, \varepsilon_6-\varepsilon_5,
\varepsilon_9-\varepsilon_7, \varepsilon_9-\varepsilon_8,
\varepsilon_{10}-\varepsilon_7, \varepsilon_{10}-\varepsilon_8 \} 
$$
then the book of shapes is
$$
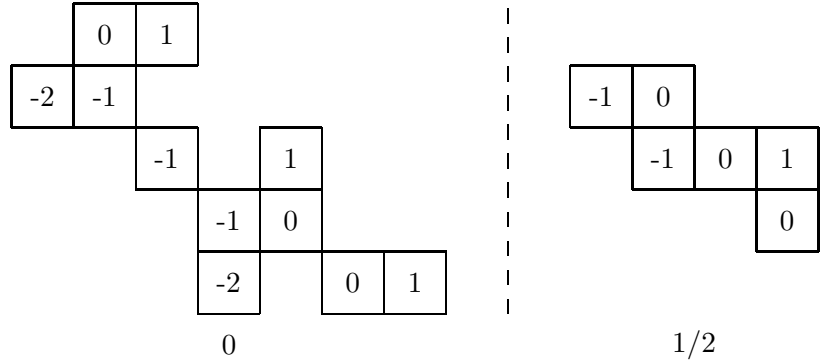
\beginpicture
\setcoordinatesystem units <0.75cm,0.75cm>         % sets scale
\setplotarea x from 0 to 13, y from -1 to 5   % sets plot size up
\linethickness=0.5pt                          % sets line thickness
\multiput{-1} at 1.5 3.5 *2 1 -1 /
\put{0}  at 1.5 4.5    %
\put{0}  at 4.5 1.5    %
\put{0}  at 5.5 0.5    %  puts contents in first page
\put{-2} at 0.5 3.5    %
\put{-2} at 3.5 0.5    %
\put{1}  at 2.5 4.5    %
\put{1}  at 4.5 2.5    %
\put{1}  at 6.5 0.5    %
\putrule from 3 0 to 4 0          %            
\putrule from 5 0 to 7 0          %  draws horizontal lines         
\putrule from 3 1 to 7 1          %   on first page        
\putrule from 2 2 to 5 2          %           
\putrule from 0 3 to 3 3          %             
\putrule from 4 3 to 5 3          %           
\putrule from 0 4 to 3 4          %           
\putrule from 1 5 to 3 5          %             
\putrule from 0 3 to 0 4        %             
\putrule from 1 3 to 1 5        %  draws vertical lines           
\putrule from 2 2 to 2 5        %   on first page           
\putrule from 3 0 to 3 3        %  
\putrule from 3 4 to 3 5        %  
\putrule from 4 0 to 4 3        %  
\putrule from 5 0 to 5 3        %  
\putrule from 6 0 to 6 1        %  
\putrule from 7 0 to 7 1        %  
\put{0}  at 10.5 3.5    %
\put{0}  at 11.5 2.5    %
\put{0}  at 12.5 1.5    %  puts contents in second page
\put{-1} at 9.5  3.5    %
\put{-1} at 10.5 2.5    %
\put{1}  at 12.5 2.5    %
\putrule from 12 1 to 13 1          %            
\putrule from 10 2 to 13 2          %  draws horizontal lines         
\putrule from 9  3 to 13 3          %   on second page        
\putrule from 9  4 to 11 4          %           
\putrule from 9  3 to 9  4        %             
\putrule from 10 2 to 10 4        %  draws vertical lines           
\putrule from 11 2 to 11 4        %   on second page           
\putrule from 12 1 to 12 3        %  
\putrule from 13 1 to 13 3        %  
\setdashes
\putrule from 8 0 to 8 5           % draws break between pages
\put{0}   at 3.5 -0.5      %
\put{1/2} at 11 -0.5       % labels pages
\endpicture
$$
where the numbers in the boxes are the contents of the boxes.
The filling 
$$
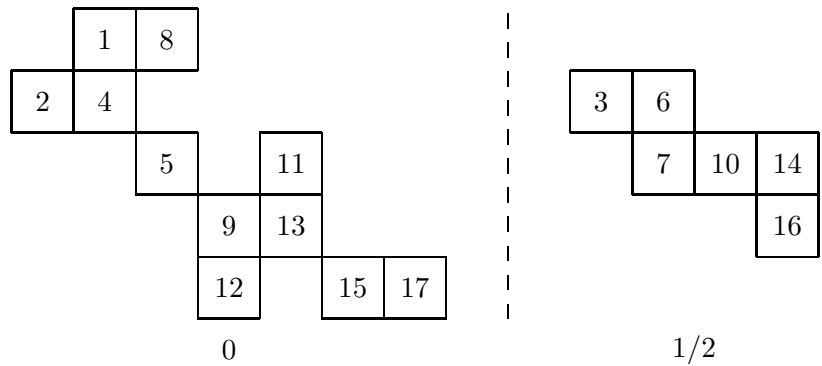
\beginpicture
\setcoordinatesystem units <0.75cm,0.75cm>         % sets scale
\setplotarea x from 0 to 13, y from -1 to 5   % sets plot size up
\linethickness=0.5pt                          % sets line thickness
\put{4}  at 1.5 3.5    %
\put{5}  at 2.5 2.5    %
\put{9}  at 3.5 1.5    %
\put{2}  at 0.5 3.5    %
\put{12} at 3.5 0.5    %
\put{1}  at 1.5 4.5    %
\put{13} at 4.5 1.5    %
\put{15} at 5.5 0.5    %  puts fillings in first page
\put{8}  at 2.5 4.5    %
\put{11} at 4.5 2.5    %
\put{17} at 6.5 0.5    %
\putrule from 3 0 to 4 0          %            
\putrule from 5 0 to 7 0          %  draws horizontal lines         
\putrule from 3 1 to 7 1          %   on first page        
\putrule from 2 2 to 5 2          %           
\putrule from 0 3 to 3 3          %             
\putrule from 4 3 to 5 3          %           
\putrule from 0 4 to 3 4          %           
\putrule from 1 5 to 3 5          %             
\putrule from 0 3 to 0 4        %             
\putrule from 1 3 to 1 5        %  draws vertical lines           
\putrule from 2 2 to 2 5        %   on first page           
\putrule from 3 0 to 3 3        %  
\putrule from 3 4 to 3 5        %  
\putrule from 4 0 to 4 3        %  
\putrule from 5 0 to 5 3        %  
\putrule from 6 0 to 6 1        %  
\putrule from 7 0 to 7 1        %  
\put{6}  at 10.5 3.5    %
\put{10} at 11.5 2.5    %
\put{16} at 12.5 1.5    %  puts fillings in second page
\put{3}  at 9.5  3.5    %
\put{7}  at 10.5 2.5    %
\put{14} at 12.5 2.5    %
\putrule from 12 1 to 13 1          %            
\putrule from 10 2 to 13 2          %  draws horizontal lines         
\putrule from 9  3 to 13 3          %   on second page        
\putrule from 9  4 to 11 4          %           
\putrule from 9  3 to 9  4        %             
\putrule from 10 2 to 10 4        %  draws vertical lines           
\putrule from 11 2 to 11 4        %   on second page           
\putrule from 12 1 to 12 3        %  
\putrule from 13 1 to 13 3        %  
\setdashes
\putrule from 8 0 to 8 5           % draws break between pages
\put{0}   at 3.5 -0.5      %
\put{1/2} at 11 -0.5       % labels pages
\endpicture
$$
is a standard tableau of shape $(\vec{\gamma},J)$.  This filling corresponds to
the permutation 
$$w=(2,12,4,5,9,1,13,15,8,11,17,3,7,6,10,16,14)
\qquad\hbox{in $\cF^{(\vec\gamma,J)}\subseteq S_{16}$.}
\qquad\qquad\hbox{\qed}$$

\section 4.  Skew shapes, ribbons, conjugation, etc. in Type A

As in the previous section let $R$ be the root system of Type
$A_{n-1}$ as given in (3.1).
For clarity, we shall state all of the results in this section
for placed shapes $(\gamma,J)$ such that $\gamma$ is dominant
and integral, i.e. $\gamma=(\gamma_1,\ldots,\gamma_n)$
with $\gamma_1\le \cdots\gamma_n$ and $\gamma_i\in \ZZ$.  
This simplification is mathematically unimportant, the reason for it
is that it allows us to avoid the notational difficulties
which arise when one wants to use books of placed shapes with several pages.

\subsection{Axial distance}

Let $(\gamma,J)$ be a placed shape such that $\gamma$ is dominant
and integral.  Let $w\in \cF^{(\gamma,J)}$ and let $t$ be the corresponding
standard tableau as defined by the map in Theorem 3.5.  
Then it follows from the definitions of $\gamma$ and $w_t$ in 
(3.2) and (3.3) that
$$\langle w\gamma,\varepsilon_i\rangle
=\langle \gamma, w^{-1}\varepsilon_i\rangle
=c({\rm box}_{w^{-1}(i)})=c(t(i)),
$$
where $t(i)$ is the box of $t$ containing the entry $i$.

In classical standard tableau theory the {\it axial distance} between
two boxes in a standard tableau is defined as follows.
Let $\lambda$ be a partition and let $t$ be a standard tableau of 
shape $\lambda$.  Let $1\le i,j\le n$ and let $t(i)$ and $t(j)$
be the boxes which are filled with $i$ and $j$ respectively.
Let $(r_i,c_i)$ and $(r_j,c_j)$ be the positions of these boxes,
where the rows and columns of $\lambda$ are numbered in the same way as for 
matrices.  Then the {\it axial distance} from $j$ to $i$ in $t$ is
$$d_{ji}(t) = c_j-c_i+r_i-r_j,$$
(see [Wz]).  This may seem confusing at first
but it is simpler if we rewrite it in terms of the corresponding 
placed shape $(\gamma,J)$ where $\gamma$ is the sequence in 
$\RR^n$ determined by some placing of $\lambda$ on infinite graph paper.
Then one gets 
$$d_{ji}(t) = c(t(j))-c(t(i)) 
=\langle w\gamma, \varepsilon_j-\varepsilon_i\rangle
=d_{\varepsilon_j-\varepsilon_i}(w),$$
where $w\in \cF^{(\gamma,J)}$ is the permutation corresponding to the standard
tableau $t$ (see (3.2)) and $d_\alpha(w)$ is the generalized axial distance
defined in (1.10).  This shows that the axial distance defined in
(1.10) is a generalization of the classical notion of axial distance.

\subsection{Skew shapes.}

The following proposition shows that, in the case of a root system
of type $A$, the definition of generalized skew shape coincides with 
the classical notion of a skew shape.  

\prop  
Let $(\gamma,J)$ be a placed shape such that $\gamma$ is 
dominant and integral. Then
the configuration of boxes associated to $(\gamma,J)$ is a placed
skew shape if and only if $(\gamma,J)$ is a skew shape as defined in (1.4).
\endprop
\pf
$\Longleftarrow:$  We shall show that if the placed configuration corresponding
to the pair $(\gamma,J)$ has any $2\times 2$ blocks of the forms
$$
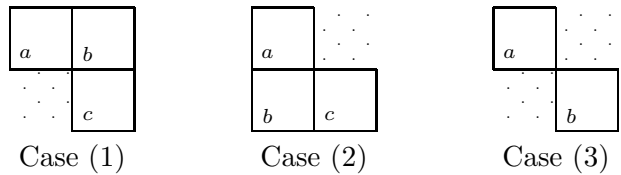
\matrix{
\beginpicture
\setcoordinatesystem units <0.75cm,0.75cm>         % sets scale
\setplotarea x from 0 to 2, y from 0 to 2    % sets plot size up
\linethickness=0.5pt                          % sets line thickness
\put{$\scriptstyle{a}$} at 0.25 1.25
\put{$\scriptstyle{b}$}  at 1.25 1.25
\put{$\scriptstyle{c}$} at 1.25 0.25
\putrule from 0 2 to 2 2          %            
\putrule from 0 1 to 2 1          %  draws horizontal lines         
\putrule from 1 0 to 2 0          %           
\putrule from 0 1 to 0 2        %             
\putrule from 1 0 to 1 2        %             
\putrule from 2 0 to 2 2        %  draws vertical lines           
\vshade 0 0 1    1 0 1  /           
\endpicture
&\qquad
&\beginpicture
\setcoordinatesystem units <0.75cm,0.75cm>         % sets scale
\setplotarea x from 0 to 2, y from 0 to 2    % sets plot size up
\linethickness=0.5pt                          % sets line thickness
\put{$\scriptstyle{a}$} at 0.25 1.25
\put{$\scriptstyle{b}$}  at 0.25 0.25
\put{$\scriptstyle{c}$} at 1.25 0.25
\putrule from 0 2 to 1 2          %            
\putrule from 0 1 to 2 1          %  draws horizontal lines         
\putrule from 0 0 to 2 0          %           
\putrule from 0 0 to 0 2        %             
\putrule from 1 0 to 1 2        %             
\putrule from 2 0 to 2 1        %  draws vertical lines           
\vshade 1 1 2    2 1 2  /           
\endpicture
&\qquad
&\beginpicture
\setcoordinatesystem units <0.75cm,0.75cm>         % sets scale
\setplotarea x from 0 to 2, y from 0 to 2    % sets plot size up
\linethickness=0.5pt                          % sets line thickness
\put{$\scriptstyle{a}$} at 0.25 1.25
\put{$\scriptstyle{b}$}  at 1.25 0.25
\putrule from 0 2 to 1 2          %            
\putrule from 0 1 to 2 1          %  draws horizontal lines         
\putrule from 1 0 to 2 0          %           
\putrule from 0 1 to 0 2        %             
\putrule from 1 0 to 1 2        %             
\putrule from 2 0 to 2 1        %  draws vertical lines           
\vshade 0 0 1    1 0 1  /           
\vshade 1 1 2    2 1 2  /           
\endpicture
\cr
\hbox{Case (1)} 
&&\hbox{Case (2)}
&&\hbox{Case (3)} \cr
}
$$
then $w\gamma\,\big\vert_K$ is not regular for some appropriate $w\in
\cF^{(\gamma,J)}$ and subsystem $R_K$ in $R$.  This will show that the placed
configuration must be a placed skew shape if $(\gamma,J)$ is a
generalized skew shape. In the pictures above the shaded regions indicate the
absence of a box and, for notational reference, we have labeled the 
boxes with $a,b,c$.

\medskip\noindent
{\it Case} (1):  Create a standard
tableau $t$ such that the $2\times 2$ block is filled with
$$
\beginpicture
\setcoordinatesystem units <0.75cm,0.75cm>         % sets scale
\setplotarea x from 0 to 2, y from 0 to 2    % sets plot size up
\linethickness=0.5pt                          % sets line thickness
\put{$i-1$} at 0.5 1.5
\put{$i$}  at 1.5 1.5
\put{$i+1$} at 1.5 0.5
\putrule from 0 2 to 2 2          %            
\putrule from 0 1 to 2 1          %  draws horizontal lines         
\putrule from 1 0 to 2 0          %           
\putrule from 0 1 to 0 2        %             
\putrule from 1 0 to 1 2        %             
\putrule from 2 0 to 2 2        %  draws vertical lines           
\vshade 0 0 1    1 0 1  /           
\endpicture
$$
by filling the region of the configuration strictly north and weakly 
west of box c in row reading order (sequentially left to right across
the rows starting at the top), putting the next entry in box c, 
and filling the remainder of the
configuration in column reading order (sequentially down the columns
beginning at the leftmost available column).  
Let $w=w(t)$ be the permutation in $\cF^{(\gamma,J)}$
which corresponds to the standard tableau $t$.  Let $t(i)$ denote the
box containing $i$ in $t$. Then, using the first identity (4.1),
$$\langle w\gamma,\alpha_i+\alpha_{i+1}\rangle
=\langle w\gamma,\varepsilon_{i+1}-\varepsilon_{i-1}\rangle
=c(t(i+1))-c(t(i-1)) = 0,$$
since the boxes $t(i+1)$ and $t(i-1)$ are on the same diagonal.
It follows that $w\gamma\,\big\vert_{\{\alpha_{i+1},\alpha_i\}}$ is not regular.

\medskip\noindent
{\it Case} (2):   
Create a standard tableau $t$ such that the $2\times 2$ block is
filled with
$$
\beginpicture
\setcoordinatesystem units <0.75cm,0.75cm>         % sets scale
\setplotarea x from 0 to 2, y from 0 to 2    % sets plot size up
\linethickness=0.5pt                          % sets line thickness
\put{$i-1$} at 0.5 1.5
\put{$i$}  at 0.5 0.5
\put{$i+1$} at 1.5 0.5
\putrule from 0 2 to 1 2          %            
\putrule from 0 1 to 2 1          %  draws horizontal lines         
\putrule from 0 0 to 2 0          %           
\putrule from 0 0 to 0 2        %             
\putrule from 1 0 to 1 2        %             
\putrule from 2 0 to 2 1        %  draws vertical lines           
\vshade 1 1 2    2 1 2  /           
\endpicture
$$
by filling the region weakly north and strictly west of box c in column reading
order, putting the next entry in box c, and filling the remainder of the
configuration in row reading order.  Using this standard tableau $t$,
the remainder of the argument is the same as for case (1).

\medskip\noindent
{\it Case} (3): 
Create a standard tableau $t$ such that the $2\times 2$ block is
filled with
$$
\beginpicture
\setcoordinatesystem units <0.75cm,0.75cm>         % sets scale
\setplotarea x from 0 to 2, y from 0 to 3    % sets plot size up
\linethickness=0.5pt                          % sets line thickness
\put{$i-1$} at 0.5 1.5
\put{$i$}  at 1.5 0.5
\putrule from 0 2 to 1 2          %            
\putrule from 0 1 to 2 1          %  draws horizontal lines         
\putrule from 1 0 to 2 0          %           
\putrule from 0 1 to 0 2        %             
\putrule from 1 0 to 1 2        %             
\putrule from 2 0 to 2 1        %  draws vertical lines           
\vshade 0 0 1    1 0 1  /           
\vshade 1 1 2    2 1 2  /           
\endpicture
$$
by filling the region strictly north and strictly west of box b in 
column reading order, putting the next entry in box b, and 
filling the remainder of the
configuration in row reading order.  Let $w=w(t)$ be the permutation
in $\cF^{(\gamma,J)}$ corresponding to $t$ and let $t(i)$ denote the
box containing $i$ in $t$. Then 
$$\langle w\gamma,\alpha_i\rangle
=\langle w\gamma,\varepsilon_i-\varepsilon_{i-1}\rangle
=c(t(i))-c(t(i-1)) = 0,$$
since $t(i)$ and $t(i-1)$ are on the same diagonal.
It follows that $w\gamma\,\big\vert_{\{\alpha_i\}}$ is not regular.

$\Longrightarrow:$  
Let $\gamma\in \ZZ^n$ and $\lambda/\mu$ describe a placed skew
shape (a skew shape placed on infinite graph paper).
Let $(\gamma,J)$ be the corresponding (generalized) placed shape as defined in
(3.3).  Let $w\in \cF^{(\gamma,J)}$ and let $t$ be the corresponding
standard tableau of shape $\lambda/\mu$.  Consider a 
$2\times 2$ block of boxes of $t$.
If these boxes are filled with
$$
\beginpicture
\setcoordinatesystem units <0.75cm,0.75cm>         % sets scale
\setplotarea x from 0 to 4, y from 0 to 3    % sets plot size up
\linethickness=0.5pt                          % sets line thickness
\put{$i$} at 0.5 1.5
\put{$j$}  at 1.5 1.5
\put{$k$}  at 0.5 0.5
\put{$\ell$} at 1.5 0.5
\putrule from 0 2 to 2 2          %            
\putrule from 0 1 to 2 1          %  draws horizontal lines         
\putrule from 0 0 to 2 0          %           
\putrule from 0 0 to 0 2        %             
\putrule from 1 0 to 1 2        %             
\putrule from 2 0 to 2 2        %  draws vertical lines           
%\vshade 0 0 1    1 0 1  /           
\endpicture
$$
then either $i<j<k<\ell$ or $i<k<j<\ell$.  In either case we have
$i<\ell-1$ and it follows that $\ell-1$ and $\ell$ are not on the same diagonal.
Thus
$$\langle w\gamma,\alpha_\ell\rangle=c(t(\ell))-c(t(\ell-1))\ne 0,$$
and so $w\gamma\,\big\vert_{\{\alpha_\ell\}}$ is regular.

The same argument shows that one can never get a standard tableau in which
$\ell$ and $\ell-2$ occur in adjacent boxes of the same diagonal and thus
it follows that $w\gamma\,\big\vert_{\{\alpha_{\ell-1},\alpha_\ell\}}$ is 
regular for all $w\in \cF^{(\gamma,J)}$.

Thus $(\gamma,J)$ is a placed skew shape in the sense of (1.4).
\endpf

\subsection{Ribbon Shapes.}

In classical tableaux theory a {\it border strip} (or {\it ribbon}) is a skew
shape which contains at most one box in each diagonal.  Although the convention,
[Mac, I \S 1 p. 5], is to assume that border strips are connected skew
shapes we shall {\it not} assume this.

Recall from (1.5) that a placed shape $(\gamma,J)$ is a
placed {\it ribbon} shape if $\gamma$ is regular, 
i.e. $\langle \gamma,\alpha\rangle \ne 0$ for all $\alpha\in R$.

\prop 
Let $(\gamma,J)$ be a placed ribbon shape such that $\gamma$
is dominant and integral. 
Then the configuration of boxes coresponding
to $(\gamma,J)$ is a placed border strip.
\pf
Let $(\gamma,J)$ be a placed ribbon shape with $\gamma$ dominant
and regular.  Since $\gamma=(\gamma_1,\ldots,\gamma_n)$ is regular,
$\gamma_i\ne \gamma_j$ for all $i\ne j$.  In terms of the placed 
configuration $\gamma_i=c({\rm box}_i)$ is the diagonal that ${\rm box}_i$ 
is on.  Thus the configuration of boxes corresponding to $(\gamma,J)$ contains
at most one box in each diagonal.
\endpf

\noindent
{\sl Example.}
If $\gamma=(-6,-5,-4,0,1,3,4,5,6,7)$ and
$J=\{\varepsilon_2-\varepsilon_1,
\varepsilon_5-\varepsilon_4,
\varepsilon_7-\varepsilon_6,
\varepsilon_9-\varepsilon_8,
\varepsilon_{10}-\varepsilon_9\}$
then the placed configuration of boxes corresponding to $(\gamma,J)$
is the placed border strip
$$
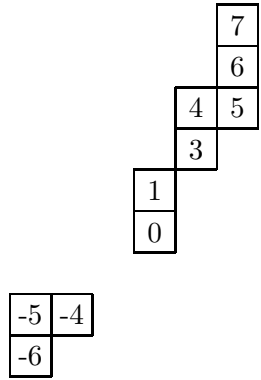
\beginpicture
\setcoordinatesystem units <0.5cm,0.5cm> point at 6 0       % sets scale
\setplotarea x from 0 to 6, y from -1 to 9    % sets plot size up
\linethickness=0.5pt                          % sets line thickness
%\put{contents of boxes}  at 2.5 -0.5    % 
\put{-6}  at 0.5 0.5    % 
\put{-5}  at 0.5 1.5    % 
\put{-4}  at 1.5 1.5    % 
\put{0}  at 3.5 3.5    % 
\put{1}  at 3.5 4.5    %  labels contents of boxes
\put{3}  at 4.5 5.5    % 
\put{4}  at 4.5 6.5    % 
\put{5}  at 5.5 6.5    % 
\put{6}  at 5.5 7.5    % 
\put{7}  at 5.5 8.5    % 
\putrule from 0 0 to 1 0          %            
\putrule from 0 1 to 2 1          %            
\putrule from 0 2 to 2 2          %            
\putrule from 3 3 to 4 3          %            
\putrule from 3 4 to 4 4          %            
\putrule from 3 5 to 5 5          %        
\putrule from 4 6 to 6 6          %           
\putrule from 4 7 to 6 7          %           
\putrule from 5 8 to 6 8          %           
\putrule from 5 9 to 6 9          %           
\putrule from 0 0 to 0 2        %             
\putrule from 1 0 to 1 2        %             
\putrule from 2 1 to 2 2        %             
\putrule from 3 3 to 3 5        %  draws vertical lines           
\putrule from 4 3 to 4 7        %             
\putrule from 5 5 to 5 9        %  
\putrule from 6 6 to 6 9        %  
\endpicture
$$
where we have labeled the boxes with their contents.
\endpf

\subsection{Conjugation of Shapes.}

Let $(\gamma,J)$ be a placed shape with $\gamma$ dominant and integral
(i.e. $\gamma=(\gamma_1,\ldots,\gamma_n)$ with 
$\gamma_1\leq\cdots\leq\gamma_n$ and $\gamma_i\in\ZZ$)
and view $(\gamma,J)$ as a placed configuration of boxes.  
In terms of placed configurations,
conjugation of shapes is equivalent to transposing
the placed configuration across the diagonal of boxes of content $0$.
The following example illustrates this.

\bigskip\noindent
{\sl Example.}
Suppose $\gamma=(-1,-1,-1,0,0,1,1)$ and $J=(\varepsilon_4-\varepsilon_2,
\varepsilon_4-\varepsilon_3, \varepsilon_6-\varepsilon_5,
\varepsilon_7-\varepsilon_5)$. Then the placed configuration of boxes
corresponding to  $(\gamma,J)$ is 
$$
\beginpicture
\setcoordinatesystem units <0.75cm,0.75cm>         % sets scale
\setplotarea x from 0 to 4, y from 0 to 3    % sets plot size up
\linethickness=0.5pt                          % sets line thickness
\put{-1} at 0.5 2.5
\put{0}  at 1.5 2.5
\put{-1} at 1.5 1.5
\put{-1} at 2.5 0.5
\put{1}  at 2.5 2.5
\put{1}  at 3.5 1.5
\put{0}  at 3.5 0.5
\putrule from 0 3 to 3 3          %            
\putrule from 0 2 to 4 2          %  draws horizontal lines         
\putrule from 1 1 to 4 1          %           
\putrule from 2 0 to 4 0          %           
\putrule from 0 2 to 0 3        %             
\putrule from 1 1 to 1 3        %             
\putrule from 2 0 to 2 3        %  draws vertical lines           
\putrule from 3 0 to 3 3        %             
\putrule from 4 0 to 4 2        %  
\vshade 2 1 2   3 1 2 /           
\endpicture
$$
in which the shaded box is not a box in the configuration.

The minimal length representative of the coset $w_0W_\gamma$ is
the permutation
$$u=\left(
\matrix{
1&2&3&4&5&6&7\cr
5&6&7&3&4&1&2 \cr
}
\right).$$ 
We have $-u\gamma=-w_0\gamma=(-1,-1,0,0,1,1,1)$ and 
$$
\eqalign{
-u(P(\gamma)\setminus J) &= 
-u\left\{ \varepsilon_4-\varepsilon_1,
\varepsilon_5-\varepsilon_1, \varepsilon_5-\varepsilon_2,
\varepsilon_5-\varepsilon_3, \varepsilon_6-\varepsilon_4,
\varepsilon_7-\varepsilon_4 \right\} \cr
&= -\left\{ \varepsilon_3-\varepsilon_5,
\varepsilon_4-\varepsilon_5, \varepsilon_4-\varepsilon_6,
\varepsilon_4-\varepsilon_7, \varepsilon_1-\varepsilon_3,
\varepsilon_2-\varepsilon_3 \right\} \cr
&= \left\{ \varepsilon_5-\varepsilon_3,
\varepsilon_5-\varepsilon_4, \varepsilon_6-\varepsilon_4,
\varepsilon_7-\varepsilon_4, \varepsilon_3-\varepsilon_1,
\varepsilon_3-\varepsilon_2 \right\}. \cr
}
$$
Thus the configuration of boxes corresponding to the placed shape 
$(\gamma,J)'$ is 
$$
\beginpicture
\setcoordinatesystem units <0.75cm,0.75cm>         % sets scale
\setplotarea x from 0 to 3, y from 0 to 4    % sets plot size up
\linethickness=0.5pt                          % sets line thickness
\put{0}  at 2.5 0.5 
\put{1}  at 2.5 1.5 
\put{-1} at 1.5 0.5 
\put{1}  at 1.5 2.5 
\put{-1} at 0.5 1.5 
\put{0}  at 0.5 2.5 
\put{1}  at 0.5 3.5 
\putrule from 3 0 to 3 2          %            
\putrule from 2 0 to 2 3          %  draws vertical lines         
\putrule from 1 0 to 1 4          %           
\putrule from 0 1 to 0 4          %           
\putrule from 1 0 to 3 0        %             
\putrule from 0 1 to 3 1        %             
\putrule from 0 2 to 3 2        %  draws horizontal lines           
\putrule from 0 3 to 2 3        %             
\putrule from 0 4 to 1 4        %  
\vshade 1 1 2   2 1 2 /           
\endpicture
$$
\endexample

\subsection{Row reading and column reading tableaux.}

Let $(\gamma,J)$ be a placed shape such that $\gamma$ is 
dominant and integral and consider the placed configuration of 
boxes corresponding to $(\gamma,J)$.  
The {\it minimal box} of the configuration is the box such that
\smallskip
\itemitem{($m_1$)} there is no box immediately above,
\smallskip
\itemitem{($m_2$)} there is no box immediately to the left, 
\smallskip
\itemitem{($m_3$)} there is no box northwest in the same diagonal, and
\smallskip
\itemitem{($m_4$)} it has the minimal content of the boxes satisfying 
($m_1$), ($m_2$) and ($m_3$).
\smallskip\noindent
There is at most one box in each diagonal satisfying ($m_1$), ($m_2$),
and ($m_3$).  Thus, ($m_4$) guarantees that the minimal box is unique.
It is clear that the minimal box of the configuration
always exists.

The {\it column reading} tableaux of shape $(\gamma,J)$ is
the filling $t_{\rm min}$ which is created inductively by 
\smallskip
\item{(a)} filling the minimal box of the configuration with $1$, and
\smallskip
\item{(b)} if $1,2,\ldots, i$ have been filled in then fill the minimal 
box of the configuration formed by the unfilled boxes with $i+1$.
\smallskip\noindent
The {\it row reading tableau} of shape $(\gamma,J)$ is the
standard tableau $t_{\rm max}$ whose conjugate 
$(t_{\rm max})'$ is the column reading tableaux
for the shape $(\gamma,J)'$ (the conjugate shape to $(\gamma,J)$).

Recall the definitions of the weak Bruhat order and
closed subsets of roots given in (1.11).

\thm  Let $(\gamma,J)$ be a placed shape such that $\gamma$ is 
dominant and integral
(i.e. $\gamma=(\gamma_1,\ldots,\gamma_n)$ with 
$\gamma_1\leq\cdots\leq\gamma_n$ and $\gamma_i\in\ZZ$).
Let $t_{\rm min}$ and $t_{\rm max}$
be the column reading and row reading tableaux of shape $(\gamma,J)$,
respectively, and let $w_{\rm min}$ and $w_{\rm max}$ be the corresponding
permutations in $\cF^{(\gamma,J)}$.
Then 
$$R(w_{\rm min})=\overline{J},
\qquad
R(w_{\rm max})= \overline{(P(\gamma)\setminus J)\cup Z(\gamma)}^c,
\qquad\hbox{and}\qquad
\cF^{(\gamma,J)} = 
[w_{min},w_{max}],$$
where $K^c$ denotes the complement of $K$ in $R^+$ and $[w_{min},w_{max}]$ denotes
the interval between $w_{\rm min}$ and $\,w_{\rm max}$ in the weak Bruhat order.
\pf
(a)  Consider the configuration of boxes corresponding to $(\gamma,J)$.
If $k>i$ then either
$c({\rm box}_k)>c({\rm box}_i)$, or
${\rm box}_k$ is in the same diagonal and southeast of 
${\rm box}_i$.
Thus when we create $t_{\rm min}$ we have that
$$\hbox{If $k>i$ then ${\rm box}_k$ gets filled before
${\rm box}_i$ if and only if ${\rm box}_k$ is {\sl northwest} of
${\rm box}_i$,}
$$
where the {\sl northwest} is in a very strong sense:
There is a sequence of boxes
$${\rm box}_i={\rm box}_{i_0},\enspace {\rm box}_{i_1},\enspace\ldots, 
\enspace{\rm box}_{i_r}={\rm box}_k$$
such that ${\rm box}_{i_m}$ is either directly above ${\rm box}_{i_{m-1}}$
or in the same diagonal and directly northwest of ${\rm box}_{i_{m-1}}$.
In other words,
$$
\hbox{ If $k>i$ then $t_{\rm min}({\rm box}_k)<t_{\rm min}({\rm box}_i)$
\quad$\Longleftrightarrow$\quad ${\rm box}_k$ is {\sl northwest} of ${\rm
box}_i$.}
$$
So, from the formula for $w_t$ in (3.2) we get
$$
\hbox{If $k>i$ then $w_{\rm min}(k)<w_{\rm min}(i)$
$\quad\Longleftrightarrow\quad \varepsilon_k-\varepsilon_i\in \overline{J}$, }
$$
where $w_{\rm min}$ is the permutation in $\cF^{(\gamma,J)}$
which corresponds to the filling $t_{\rm min}$  
and $\overline{J}$ is the closure of $J$ in $R$.
It follows that 
$$R(w_{\rm min})=\overline{J}.$$

(b) There are at least two ways to prove that
$R(w_{\rm max})= \overline{(P(\gamma)\setminus J)\cup Z(\gamma)}^c$.
One can mimic the proof of part (a) by defining the maximal
box of a configuration and a corresponding filling.  
Alternatively one can use the definition of conjugation 
and the fact that $R(w_0w)=R(w)^c$.
The permutation $w_{\rm min}$ is the unique minimal element
of $\cF^{(\gamma,J)}$ and the conjugate of $w_{\rm max}$ is the
unique minimal element of $\cF^{(\gamma,J)'}$.
We shall leave the details to the reader.

(c)  An element $w\in W$ 
is an element of $\cF^{(\gamma,J)}$ if and only if
$R(w)\cap P(\gamma)=J$ and $R(w)\cap Z(\gamma)=\emptyset$.
Thus $\cF^{(\gamma,J)}$ consists of those permutations
$w\in W$ such that 
$$\overline{J}\subseteq R(w) \subseteq 
\overline{(P(\gamma)\setminus J)\cup Z(\gamma)}^c.$$
Since the weak Bruhat order is the ordering determined
by inclusions of $R(w)$, it follows that $\cF^{(\gamma,J)}$ is
the interval between $w_{\rm min}$ and $w_{\rm max}$.
\endpf

\noindent
{\sl Example.}
Suppose $\gamma=(-1,-1,-1,0,0,1,1)$ and $J=\{\varepsilon_4-\varepsilon_2,
\varepsilon_4-\varepsilon_3, \varepsilon_6-\varepsilon_5,
\varepsilon_7-\varepsilon_5\}$. The minimal and maximal elements in 
$\cF^{(\gamma,J)}$ are the 
permutations
$$
w_{\rm min}=\pmatrix{
1&2&3&4&5&6&7\cr
1&3&4&2&7&5&6\cr}
\qquad\hbox{and}\qquad
w_{\rm max}=\pmatrix{
1&2&3&4&5&6&7\cr
1&5&6&2&7&3&4\cr}.
$$
The permutations correspond to the standard tableaux
$$
\beginpicture
\setcoordinatesystem units <0.75cm,0.75cm>         % sets scale
\setplotarea x from 0 to 4, y from 0 to 3    % sets plot size up
\linethickness=0.5pt                          % sets line thickness
\put{1} at 0.5 2.5
\put{2} at 1.5 2.5
\put{3} at 1.5 1.5
\put{4} at 2.5 0.5
\put{5} at 2.5 2.5
\put{6} at 3.5 1.5
\put{7} at 3.5 0.5
\putrule from 0 3 to 3 3          %            
\putrule from 0 2 to 4 2          %  draws horizontal lines         
\putrule from 1 1 to 4 1          %           
\putrule from 2 0 to 4 0          %           
\putrule from 0 2 to 0 3        %             
\putrule from 1 1 to 1 3        %             
\putrule from 2 0 to 2 3        %  draws vertical lines           
\putrule from 3 0 to 3 3        %             
\putrule from 4 0 to 4 2        %  
\vshade 2 1 2   3 1 2 /           
\endpicture
\qquad\qquad\hbox{and}\qquad\qquad
\beginpicture
\setcoordinatesystem units <0.75cm,0.75cm>         % sets scale
\setplotarea x from 0 to 4, y from 0 to 3    % sets plot size up
\linethickness=0.5pt                          % sets line thickness
\put{1} at 0.5 2.5
\put{2} at 1.5 2.5
\put{5} at 1.5 1.5
\put{6} at 2.5 0.5
\put{3} at 2.5 2.5
\put{4} at 3.5 1.5
\put{7} at 3.5 0.5
\putrule from 0 3 to 3 3          %            
\putrule from 0 2 to 4 2          %  draws horizontal lines         
\putrule from 1 1 to 4 1          %           
\putrule from 2 0 to 4 0          %           
\putrule from 0 2 to 0 3        %             
\putrule from 1 1 to 1 3        %             
\putrule from 2 0 to 2 3        %  draws vertical lines           
\putrule from 3 0 to 3 3        %             
\putrule from 4 0 to 4 2        %  
\vshade 2 1 2   3 1 2 /           
\endpicture
.
$$
\endexample

\section 5. Standard tableaux for type $C$ in terms of boxes

\subsection{The root system.}

Let $\{\varepsilon_1,\ldots,\varepsilon_n\}$ be an orthonormal basis  of
$\RR^n$ so that each sequence $\gamma=(\gamma_1,\ldots,\gamma_n)\in \RR^n$
is identified with the vector $\gamma=\sum_i \gamma_i\varepsilon_i$.
The root system of type $C_n$ is given by the sets
$$
R = \left\{\pm2\varepsilon_i,\varepsilon_j\pm\varepsilon_i
\ |\  1\leq i,j\leq n\right\}
\quad\hbox{and}\quad
R^+=\left\{2\varepsilon_i,\varepsilon_j-\varepsilon_i
\ |\  1\leq i<j\leq n\right\}.
$$
The simple roots are given by $\alpha_1=2\varepsilon_1$,
$\alpha_i=\varepsilon_i-\varepsilon_{i-1}$, $2\le i\le n$.
The Weyl group $W=WC_n$ is the {\it hyperoctahedral group}
of permutations of $-n,\ldots,-1,1,\ldots,n$ such that
$w(-i)=-w(i)$.  This groups acts on the $\varepsilon_i$ by the
rule $w\varepsilon_i=\varepsilon_{w(i)}$, with the convention
that $\varepsilon_{-i}=-\varepsilon_i$.

\subsection{Rearranging $\gamma$.}

The analysis in this case is analogous to the method that was used
in (3.7) to create books of placed configurations in the type A case. 
For clarity, we recommend that the reader compare the machinations below with the
case done in (3.7).

Let $\gamma\in \RR^n$.  By applying an element of the Weyl group to $\gamma$
we can rearrange the entries of $\gamma$ in increasing order
$(0\le \gamma_1\le \gamma_2\le \cdots\le \gamma_n)$.
Then, if $\gamma_i\in (x,x+1/2)$ for some integer $x$, replace $\gamma_i$ with
$-\gamma_i$.  
Next group the elements of $\gamma$ in terms of their $\ZZ$-cosets
and rearrange each group to be in increasing order.  
There are three kinds of groups which can occur:
$$
\matrix{
\hfil \beta+C_\beta&=& (\beta+z_1,\beta+z_2,\ldots,\beta+z_k),\hfil 
&\quad &\hbox{with $\beta\in (x+1/2,x+1)$ for
some $x\in \ZZ$,}\hfill \cr
&&&&\hbox{\qquad\qquad $z_i\in \ZZ_{\ge 0}$ and $z_1\le \ldots\le z_k$,}\hfill \cr
\hfil C_0 &=& (z_1\le z_2\le\cdots\le z_k), \hfil
&&\hbox{with $z_i\in \ZZ_{\ge 0}$ and $z_1\le \ldots\le z_k$,}\hfill \cr
\hfil 1/2+C_{1/2}&=& (1/2+z_1\le 1/2+z_2\le \cdots\le 1/2+z_k), \hfil
&&\hbox{with $z_i\in \ZZ_{\ge 0}$ and $z_1\le \ldots\le z_k$.}\hfill \cr
}
$$
Choose some ordering on the groups and let
$$\vec \gamma = (\beta_1+C_{\beta_1},\ldots,\beta_r+C_{\beta_r}).$$
Because these changes are obtained by applying elements of the
Weyl group, the calibration graphs corresponding to $\vec \gamma$ and
the $\gamma$ are the same.  
Thus it is sufficient to study the standard tableaux corresponding to 
$\vec\gamma$.

\subsection{Books of placed configurations.}

As in the type $A$ case, we have set partitions
$$
Z(\vec{\gamma})=\bigcup_{\beta_i}Z_{\beta_i}, \quad\quad
P(\vec{\gamma})=\bigcup_{\beta_i}P_{\beta_i}, \quad\hbox{
and, }\quad
$$
for any $J\subseteq P(\vec{\gamma})$,
$$J=\bigcup_{\beta_i}J_{\beta_i},
\quad\hbox{where $J_{\beta_i}=J\cap P_{\beta_i}$.}
$$
Each pair $(\beta+C_\beta,J_{\beta})$ is a placed shape and we may
identify $(\vec{\gamma},J)$ with the {\it book of placed shapes}
$$\left( (\beta_1+C_{\beta_1},J_{\beta_1}),\ldots,
(\beta_r+C_{\beta_r},J_{\beta_r}) \right).$$
We think of this as a book with pages numbered by the values
$\beta_1,\ldots,\beta_r$ and with a placed configuration determined by
$(\beta_i+C_{\beta_i},J_{\beta_i})$ on page $\beta_i$.
One determines the placed configurations as follows.

\subsection{Page $\beta$, $\beta\ne {1\over2}, 0$:}  

As in the type $A$ case we will place boxes on a page of infinite
graph paper which has the diagonals numbered consecutively with the elements 
of $\ZZ$, from bottom left to top right.  For each $1\le i\le n$, 
place ${\rm box}_i$ on diagonal $\vec \gamma_i-\beta$. 
The boxes on each diagonal are arranged in increasing order from
top left to bottom right.  Then
$$
\eqalign{
P_\beta &= \left\{
\varepsilon_j-\varepsilon_i \ |\ 
\hbox{$j>i$ and ${\rm box}_i$ and ${\rm box}_j$ are in adjacent diagonals} 
\right\}
\quad\hbox{and} \cr
Z_\beta &= \left\{
\varepsilon_j-\varepsilon_i \ |\  
\hbox{$j>i$ and ${\rm box}_i$ and ${\rm box}_j$ are in the same diagonal}
\right\}. \cr
}
$$
If $J_\beta\subseteq P_\beta$ arrange the boxes on adjacent diagonals 
according to the rules
\smallskip
\itemitem{($a_\beta$)}
if $\varepsilon_j-\varepsilon_i\in J_r$ place 
${\rm box}_j$ northwest of ${\rm box}_i$, and
\smallskip
\itemitem{($a'_\beta$)}
if $\varepsilon_j-\varepsilon_i\in P_r\backslash J_r$ place 
${\rm box}_j$ southeast of ${\rm box}_i$.
\smallskip\noindent
A {\it standard tableau} $t$ is a filling of the boxes with 
distinct entries from the set $\{-n,\ldots,-1,1,\ldots,n\}$
such that if $i$ appears then $-i$ does not appear and
\item{(a)} if $j>i$ and 
${\rm box}_j$ and ${\rm box}_i$ are in the same diagonal
then $t({\rm box}_i)<t({\rm box}_j)$,
\item{(b)} if $j>i$,  
${\rm box}_i$ and ${\rm box}_j$ are in adjacent diagonals
and ${\rm box}_j$ is northwest of ${\rm box}_i$ 
then $t({\rm box}_i)>t({\rm box}_j)$,
\item{(c)} if $j>i$, 
${\rm box}_i$ and ${\rm box}_j$ are in adjacent diagonals
and ${\rm box}_j$ is southeast of ${\rm box}_i$ 
then $t({\rm box}_i)<t({\rm box}_j)$.

\bigskip\noindent
{\sl Example.}
Suppose $\beta+C_\beta=\beta+(0,0,0,1,1,1,2,2,2,3,3,3)$ and
$$\eqalign{
J_\beta &=\left\{ \varepsilon_4-\varepsilon_1, \varepsilon_4-\varepsilon_2,
\varepsilon_4-\varepsilon_3, \varepsilon_5-\varepsilon_2,
\varepsilon_5-\varepsilon_3, \varepsilon_7-\varepsilon_5,
\varepsilon_7-\varepsilon_6, \varepsilon_8-\varepsilon_6,
\varepsilon_{10}-\varepsilon_9, \varepsilon_{10}-\varepsilon_8,\right. \cr
&\phantom{J = } \left.\varepsilon_{10}-\varepsilon_7,
\varepsilon_{11}-\varepsilon_9,
\varepsilon_{11}-\varepsilon_8, \varepsilon_{11}-\varepsilon_7,
\varepsilon_{12}-\varepsilon_9 \right\}. \cr
}
$$
The placed configuration of boxes corresponding to $(\beta+C_\beta,J_\beta)$
is as given below.
$$
\matrix{
\beginpicture
\setcoordinatesystem units <0.5cm,0.5cm> point at 6 0       % sets scale
\setplotarea x from 0 to 5, y from -1 to 6    % sets plot size up
\linethickness=0.5pt                          % sets line thickness
%\put{contents of boxes}  at 2.5 -0.5    % 
\put{0}  at 0.5 2.5    % 
\put{0}  at 1.5 1.5    % 
\put{0}  at 2.5 0.5    % 
\put{1}  at 0.5 3.5    % 
\put{1}  at 1.5 2.5    %  labels contents of boxes
\put{1}  at 3.5 0.5    % 
\put{2}  at 1.5 3.5    % 
\put{2}  at 3.5 1.5    % 
\put{2}  at 4.5 0.5    % 
\put{3}  at 0.5 5.5    % 
\put{3}  at 1.5 4.5    % 
\put{3}  at 4.5 1.5    % 
\putrule from 0 6 to 1 6          %            
\putrule from 0 5 to 2 5          %            
\putrule from 0 4 to 2 4          %            
\putrule from 0 3 to 2 3          %            
\putrule from 3 2 to 5 2          %        
\putrule from 0 2 to 2 2          %  draws horizontal lines         
\putrule from 1 1 to 5 1          %           
\putrule from 2 0 to 5 0          %           
\putrule from 0 2 to 0 4        %             
\putrule from 0 5 to 0 6        %             
\putrule from 1 1 to 1 6        %             
\putrule from 2 0 to 2 5        %  draws vertical lines           
\putrule from 3 0 to 3 2        %             
\putrule from 4 0 to 4 2        %  
\putrule from 5 0 to 5 2        %  
\endpicture
&\qquad
&\beginpicture
\setcoordinatesystem units <0.5cm,0.5cm> point at -1 0       % sets scale
\setplotarea x from 0 to 5, y from -1 to 6    % sets plot size up
\linethickness=0.5pt                          % sets line thickness
%\put{numbering of boxes}  at 2.5 -0.5    % 
\put{1}  at 0.5 2.5    % 
\put{2}  at 1.5 1.5    % 
\put{3}  at 2.5 0.5    % 
\put{4}  at 0.5 3.5    % 
\put{5}  at 1.5 2.5    %  labels numbers of boxes
\put{6}  at 3.5 0.5    % 
\put{7}  at 1.5 3.5    % 
\put{8}  at 3.5 1.5    % 
\put{9}  at 4.5 0.5    % 
\put{10} at 0.5 5.5    % 
\put{11} at 1.5 4.5    % 
\put{12} at 4.5 1.5    % 
\putrule from 0 6 to 1 6          %            
\putrule from 0 5 to 2 5          %            
\putrule from 0 4 to 2 4          %            
\putrule from 0 3 to 2 3          %            
\putrule from 3 2 to 5 2          %        
\putrule from 0 2 to 2 2          %  draws horizontal lines         
\putrule from 1 1 to 5 1          %           
\putrule from 2 0 to 5 0          %           
\putrule from 0 2 to 0 4        %             
\putrule from 0 5 to 0 6        %             
\putrule from 1 1 to 1 6        %             
\putrule from 2 0 to 2 5        %  draws vertical lines           
\putrule from 3 0 to 3 2        %             
\putrule from 4 0 to 4 2        %  
\putrule from 5 0 to 5 2        %  
\endpicture 
&\qquad
&\beginpicture
\setcoordinatesystem units <0.5cm,0.5cm> point at -1 0       % sets scale
\setplotarea x from 0 to 5, y from -1 to 6    % sets plot size up
\linethickness=0.5pt                          % sets line thickness
%\put{numbering of boxes}  at 2.5 -0.5    % 
\put{-9}  at 0.5 2.5    % 
\put{-7}  at 1.5 1.5    % 
\put{-5}  at 2.5 0.5    % 
\put{-12}  at 0.5 3.5    % 
\put{-8}  at 1.5 2.5    %  labels numbers of boxes
\put{3}  at 3.5 0.5    % 
\put{-2}  at 1.5 3.5    % 
\put{1}  at 3.5 1.5    % 
\put{4}  at 4.5 0.5    % 
\put{-11} at 0.5 5.5    % 
\put{-10} at 1.5 4.5    % 
\put{6} at 4.5 1.5    % 
\putrule from 0 6 to 1 6          %            
\putrule from 0 5 to 2 5          %            
\putrule from 0 4 to 2 4          %            
\putrule from 0 3 to 2 3          %            
\putrule from 3 2 to 5 2          %        
\putrule from 0 2 to 2 2          %  draws horizontal lines         
\putrule from 1 1 to 5 1          %           
\putrule from 2 0 to 5 0          %           
\putrule from 0 2 to 0 4        %             
\putrule from 0 5 to 0 6        %             
\putrule from 1 1 to 1 6        %             
\putrule from 2 0 to 2 5        %  draws vertical lines           
\putrule from 3 0 to 3 2        %             
\putrule from 4 0 to 4 2        %  
\putrule from 5 0 to 5 2        %  
\endpicture 
\cr
\hbox{contents of boxes}
&&\hbox{numbering of boxes}
&&\hbox{a standard tableau} \cr
}
$$
\endexample

\subsection{Page ${1\over2}$:}  

We will place boxes on a page of infinite
graph paper which has the diagonals numbered consecutively with the elements 
of ${1\over2}+\ZZ$, from bottom left to top right.  
For each $1\le i\le n$ such that $\vec \gamma_i\in {1\over2}+C_{1\over2}$, 
place ${\rm box}_i$ on diagonal $\vec\gamma_i$ 
and ${\rm box}_{-i}$ on diagonal $-\vec\gamma_i$.
The boxes on each diagonal are arranged in increasing order from
top left to bottom right.  With this placing of boxes we have
$$P_{1\over2} = \left\{
\matrix{
\varepsilon_j-\varepsilon_i, &\hbox{such that
$j>i$ and ${\rm box}_i$ and ${\rm box}_j$ are in adjacent diagonals} \hfill\cr
\varepsilon_j+\varepsilon_i, &\hbox{such that
${\rm box}_j$ and ${\rm box}_i$ are both in diagonal $1/2$}\hfill\cr
2\varepsilon_i, &\hbox{such that
${\rm box}_i$ is in diagonal $1/2$}\hfill\cr
}
\right\},
$$
$$Z_{1\over2} = \left\{
\matrix{
\varepsilon_j-\varepsilon_i, &\hbox{such that
$j>i$ and ${\rm box}_i$ and ${\rm box}_j$ are in the same diagonal} \hfill\cr
}
\right\}.
$$
If $J_{1\over2}\subseteq P_{1\over2}$ arrange the boxes on adjacent diagonals 
according to the rules:
\smallskip
\itemitem{($a_{1\over2}$)}
If $\varepsilon_j-\varepsilon_i\in J_{1\over2}$ place 
${\rm box}_j$ northwest of ${\rm box}_i$ and 
${\rm box}_{-i}$ northwest of ${\rm box}_{-j}$.
\smallskip
\itemitem{($a'_{1\over2}$)}
If $\varepsilon_j-\varepsilon_i\in P_{1\over2}\backslash J_{1\over2}$ place 
${\rm box}_j$ southeast of ${\rm box}_i$ and 
${\rm box}_{-i}$ southeast of ${\rm box}_{-j}$.
\smallskip
\itemitem{($b_{1\over2}$)}
If $\varepsilon_j+\varepsilon_i\in J_{1\over2}$ ($j\ge i$) place 
${\rm box}_j$ northwest of ${\rm box}_{-i}$ and 
${\rm box}_i$ northwest of ${\rm box}_{-j}$. 
\smallskip
\itemitem{($b'_{1\over2}$)}
If $\varepsilon_j+\varepsilon_i\in P_{1\over2}\backslash J_{1\over2}$ 
($j\ge i$) place 
${\rm box}_j$ southeast of ${\rm box}_{-i}$ and 
${\rm box}_i$ southeast of ${\rm box}_{-j}$. 

\smallskip\noindent
A {\it standard tableau} $t$ is a filling of the boxes with 
distinct entries from the set $\{-n,\ldots,-1,1,\ldots,n\}$
such that 
$$
t({\rm box}_i)=-t({\rm box}_{-i}) 
$$
and
\smallskip
\item{(a)} If $j>i$ and 
${\rm box}_j$ and ${\rm box}_i$ are in the same diagonal
then $t({\rm box}_i)<t({\rm box}_j)$,
\smallskip
\item{(b)} If $j>i$,  
${\rm box}_i$ and ${\rm box}_j$ are in adjacent diagonals
and ${\rm box}_j$ is northwest of ${\rm box}_i$ 
then $t({\rm box}_i)>t({\rm box}_j)$,
\smallskip
\item{(c)} If $j>i$, 
${\rm box}_i$ and ${\rm box}_j$ are in adjacent diagonals
and ${\rm box}_j$ is southeast of ${\rm box}_i$ 
then $t({\rm box}_i)<t({\rm box}_j)$.

\bigskip\noindent
{\sl Example.}  
Suppose ${1\over 2}+C_{1\over 2}=
\left({1\over2},{1\over2},{1\over2},{1\over2},
{3\over2},{3\over2},{3\over2},{3\over2},{5\over2},{5\over2},{7\over2}\right)$ 
and
$$\eqalign{
J_{1\over2} 
&=\left\{ \varepsilon_{11}-\varepsilon_{10}, \varepsilon_{10}-\varepsilon_8,
\varepsilon_9-\varepsilon_7, \varepsilon_9-\varepsilon_8,
\varepsilon_7-\varepsilon_3, \varepsilon_7-\varepsilon_4,
\varepsilon_6-\varepsilon_2, \varepsilon_6-\varepsilon_3,
\varepsilon_6-\varepsilon_4, \varepsilon_5-\varepsilon_4,\right. \cr
&\phantom{J = } \left.\varepsilon_5-\varepsilon_3,
\varepsilon_5-\varepsilon_2,
\varepsilon_1+\varepsilon_2, \varepsilon_1+\varepsilon_3,
\varepsilon_1+\varepsilon_4, 2\varepsilon_1 \right\}. \cr
}
$$
The placed configuration of boxes corresponding to 
$({1\over2}+C_{1\over2},J_{1\over2})$
is as given below.
$$
\matrix{
\beginpicture
\setcoordinatesystem units <0.5cm,0.5cm> point at 6 0       % sets scale
\setplotarea x from 0 to 8, y from -1 to 7    % sets plot size up
\linethickness=0.5pt                          % sets line thickness
%\put{contents of boxes}  at 2.5 -0.5    % 
\put{-${7\over2}$}  at 0.5 3.5    % 
\put{-${5\over2}$}  at 0.5 4.5    % 
\put{-${5\over2}$}  at 2.5 2.5    % 
\put{-${3\over2}$}  at 0.5 5.5    % 
\put{-${3\over2}$}  at 2.5 3.5    % 
\put{-${3\over2}$}  at 3.5 2.5    % 
\put{-${3\over2}$}  at 4.5 1.5    % 
\put{-${1\over2}$}  at 1.5 5.5    % 
\put{-${1\over2}$}  at 2.5 4.5    % 
\put{-${1\over2}$}  at 3.5 3.5    % 
\put{-${1\over2}$}  at 6.5 0.5    % 
\put{$1\over2$}  at 1.5 6.5    % 
\put{$1\over2$}  at 4.5 3.5    % 
\put{$1\over2$}  at 5.5 2.5    % 
\put{$1\over2$}  at 6.5 1.5    % 
\put{$3\over2$}  at 3.5 5.5    % 
\put{$3\over2$}  at 4.5 4.5    % 
\put{$3\over2$}  at 5.5 3.5    % 
\put{$3\over2$}  at 7.5 1.5    % 
\put{$5\over2$}  at 5.5 4.5    % 
\put{$5\over2$}  at 7.5 2.5    % 
\put{$7\over2$}  at 7.5 3.5    % 
\putrule from 1 7 to 2 7          %            
\putrule from 0 6 to 2 6          %            
\putrule from 3 6 to 4 6          %            
\putrule from 0 5 to 6 5          %            
\putrule from 0 4 to 1 4          %            
\putrule from 2 4 to 6 4          %            
\putrule from 7 4 to 8 4          %            
\putrule from 0 3 to 1 3          %        
\putrule from 2 3 to 6 3          %        
\putrule from 7 3 to 8 3          %        
\putrule from 2 2 to 8 2          %  draws horizontal lines         
\putrule from 4 1 to 5 1          %           
\putrule from 6 1 to 8 1          %           
\putrule from 6 0 to 7 0          %           
\putrule from 0 3 to 0 6        %             
\putrule from 1 3 to 1 7        %             
\putrule from 2 2 to 2 7        %             
\putrule from 3 2 to 3 6        %  draws vertical lines           
\putrule from 4 1 to 4 6        %             
\putrule from 5 1 to 5 5        %  
\putrule from 6 0 to 6 5        %  
\putrule from 7 0 to 7 4        %  
\putrule from 8 1 to 8 4        %  
\endpicture
&\qquad
&\beginpicture
\setcoordinatesystem units <0.5cm,0.5cm> point at 6 0       % sets scale
\setplotarea x from 0 to 8, y from -1 to 7    % sets plot size up
\linethickness=0.5pt                          % sets line thickness
%\put{contents of boxes}  at 2.5 -0.5    % 
\put{-11}  at 0.5 3.5    % 
\put{-10}  at 0.5 4.5    % 
\put{-9}  at 2.5 2.5    % 
\put{-8}  at 0.5 5.5    % 
\put{-7}  at 2.5 3.5    % 
\put{-6}  at 3.5 2.5    % 
\put{-5}  at 4.5 1.5    % 
\put{-4}  at 1.5 5.5    % 
\put{-3}  at 2.5 4.5    % 
\put{-2}  at 3.5 3.5    % 
\put{-1}  at 6.5 0.5    % 
\put{1}  at 1.5 6.5    % 
\put{2}  at 4.5 3.5    % 
\put{3}  at 5.5 2.5    % 
\put{4}  at 6.5 1.5    % 
\put{5}  at 3.5 5.5    % 
\put{6}  at 4.5 4.5    % 
\put{7}  at 5.5 3.5    % 
\put{8}  at 7.5 1.5    % 
\put{9}  at 5.5 4.5    % 
\put{10}  at 7.5 2.5    % 
\put{11}  at 7.5 3.5    % 
\putrule from 1 7 to 2 7          %            
\putrule from 0 6 to 2 6          %            
\putrule from 3 6 to 4 6          %            
\putrule from 0 5 to 6 5          %            
\putrule from 0 4 to 1 4          %            
\putrule from 2 4 to 6 4          %            
\putrule from 7 4 to 8 4          %            
\putrule from 0 3 to 1 3          %        
\putrule from 2 3 to 6 3          %        
\putrule from 7 3 to 8 3          %        
\putrule from 2 2 to 8 2          %  draws horizontal lines         
\putrule from 4 1 to 5 1          %           
\putrule from 6 1 to 8 1          %           
\putrule from 6 0 to 7 0          %           
\putrule from 0 3 to 0 6        %             
\putrule from 1 3 to 1 7        %             
\putrule from 2 2 to 2 7        %             
\putrule from 3 2 to 3 6        %  draws vertical lines           
\putrule from 4 1 to 4 6        %             
\putrule from 5 1 to 5 5        %  
\putrule from 6 0 to 6 5        %  
\putrule from 7 0 to 7 4        %  
\putrule from 8 1 to 8 4        %  
\endpicture
&\qquad
&\beginpicture
\setcoordinatesystem units <0.5cm,0.5cm> point at 6 0       % sets scale
\setplotarea x from 0 to 8, y from -1 to 7    % sets plot size up
\linethickness=0.5pt                          % sets line thickness
%\put{contents of boxes}  at 2.5 -0.5    % 
\put{-3}  at 0.5 3.5    % 
\put{-5}  at 0.5 4.5    % 
\put{8}  at 2.5 2.5    % 
\put{-7}  at 0.5 5.5    % 
\put{-2}  at 2.5 3.5    % 
\put{9}  at 3.5 2.5    % 
\put{10}  at 4.5 1.5    % 
\put{-6}  at 1.5 5.5    % 
\put{-4}  at 2.5 4.5    % 
\put{-1}  at 3.5 3.5    % 
\put{11}  at 6.5 0.5    % 
\put{-11}  at 1.5 6.5    % 
\put{1}  at 4.5 3.5    % 
\put{4}  at 5.5 2.5    % 
\put{6}  at 6.5 1.5    % 
\put{-10}  at 3.5 5.5    % 
\put{-9}  at 4.5 4.5    % 
\put{2}  at 5.5 3.5    % 
\put{7}  at 7.5 1.5    % 
\put{-8}  at 5.5 4.5    % 
\put{5}  at 7.5 2.5    % 
\put{3}  at 7.5 3.5    % 
\putrule from 1 7 to 2 7          %            
\putrule from 0 6 to 2 6          %            
\putrule from 3 6 to 4 6          %            
\putrule from 0 5 to 6 5          %            
\putrule from 0 4 to 1 4          %            
\putrule from 2 4 to 6 4          %            
\putrule from 7 4 to 8 4          %            
\putrule from 0 3 to 1 3          %        
\putrule from 2 3 to 6 3          %        
\putrule from 7 3 to 8 3          %        
\putrule from 2 2 to 8 2          %  draws horizontal lines         
\putrule from 4 1 to 5 1          %           
\putrule from 6 1 to 8 1          %           
\putrule from 6 0 to 7 0          %           
\putrule from 0 3 to 0 6        %             
\putrule from 1 3 to 1 7        %             
\putrule from 2 2 to 2 7        %             
\putrule from 3 2 to 3 6        %  draws vertical lines           
\putrule from 4 1 to 4 6        %             
\putrule from 5 1 to 5 5        %  
\putrule from 6 0 to 6 5        %  
\putrule from 7 0 to 7 4        %  
\putrule from 8 1 to 8 4        %  
\endpicture
\cr
\hbox{contents of boxes}
&&\hbox{numbering of boxes}
&&\hbox{a standard tableau} \cr
}
$$
\endexample

\subsection{Page $0$:}  

We will place boxes on a page of infinite
graph paper which has the diagonals numbered consecutively with the elements 
of $\ZZ$, from bottom left to top right.  
For each $1\le i\le n$ such that $\vec \gamma_i\in C_0$, 
\smallskip
\item{(1)} place ${\rm box}_i$ on diagonal $0$ if $\vec \gamma_i=0$, and
\smallskip
\item{(2)} if $\vec\gamma_i\ne 0$ place ${\rm box}_i$ on diagonal $\vec \gamma_i$ 
and ${\rm box}_{-i}$ on diagonal $-\vec\gamma_i$.
\smallskip\noindent
The boxes on each diagonal are arranged in increasing order from
top left to bottom right.  With this placing of boxes we have
$$P_0 = \left\{
\matrix{
\varepsilon_j-\varepsilon_i, 
&\hbox{such that
$j>i$ and ${\rm box}_i$ and ${\rm box}_j$ are in adjacent diagonals}\hfill \cr
\varepsilon_j+\varepsilon_i, 
&\hbox{such that
${\rm box}_j$ is in diagonal $1$ and ${\rm box}_i$ is in diagonal $0$}\hfill\cr
}
\right\},
$$
$$Z_0 = \left\{
\matrix{
\varepsilon_j-\varepsilon_i, &\hbox{such that
$j>i$ and ${\rm box}_i$ and ${\rm box}_j$ are in the same diagonal}\hfill \cr
\varepsilon_j+\varepsilon_i, &\hbox{such that
$j>i$ and ${\rm box}_j$ and ${\rm box}_i$ are in diagonal $0$} \hfill\cr
2\varepsilon_i, &\hbox{such that
${\rm box}_i$ is in diagonal $0$} \hfill\cr
}
\right\}.
$$
If $J_0\subseteq P_0$ then arrange the boxes on adjacent diagonals 
according to:
\smallskip
\itemitem{($a_0$)}
if $\varepsilon_j-\varepsilon_i\in J_0$ place 
${\rm box}_j$ northwest of ${\rm box}_i$ and 
if ${\rm box}_i$ is not on diagonal $0$ place
${\rm box}_{-i}$ northwest of ${\rm box}_{-j}$,
\smallskip
\itemitem{($a'_0$)}
if $\varepsilon_j-\varepsilon_i\in P_0\backslash J_0$ place 
${\rm box}_j$ southeast of ${\rm box}_i$ and 
if ${\rm box}_i$ is not on diagonal $0$ place
${\rm box}_{-i}$ southeast of ${\rm box}_{-j}$,
\smallskip
\itemitem{($b_0$)}
if $\varepsilon_j+\varepsilon_i\in J_0$ place 
${\rm box}_i$ northwest of ${\rm box}_{-j}$, 
\smallskip
\itemitem{($b'_0$)}
if $\varepsilon_j+\varepsilon_i\in P_0\backslash J_0$ place 
${\rm box}_i$ southeast of ${\rm box}_{-j}$. 
\smallskip\noindent
A {\it standard tableau} $t$ is a filling of the boxes with 
distinct entries from the set $\{-n,\ldots,-1,1,\ldots,n\}$
such that 
$$
\matrix{
t({\rm box}_i)=-t({\rm box}_{-i}), 
&\hbox{if ${\rm box}_i$ is not on the zero diagonal,}\hfil\cr
t({\rm box}_i)>0, 
&\hbox{if ${\rm box}_i$ is on the zero diagonal,}\hfil\cr
}$$
and
\smallskip
\item{(a)} if $j>i$ and 
${\rm box}_j$ and ${\rm box}_i$ are in the same diagonal
then $t({\rm box}_i)<t({\rm box}_j)$,
\smallskip
\item{(b)} if $j>i$,  
${\rm box}_i$ and ${\rm box}_j$ are in adjacent diagonals
and ${\rm box}_j$ is northwest of ${\rm box}_i$ 
then $t({\rm box}_i)>t({\rm box}_j)$,
\smallskip
\item{(c)} if $j>i$, 
${\rm box}_i$ and ${\rm box}_j$ are in adjacent diagonals
and ${\rm box}_j$ is southeast of ${\rm box}_i$ 
then $t({\rm box}_i)<t({\rm box}_j)$.

\bigskip\noindent
{\sl Example.}  
Suppose $C_0=(0,0,0,1,1,1,2)$ and
$$\eqalign{
J_0 &=\left\{ \varepsilon_4-\varepsilon_1, \varepsilon_4-\varepsilon_2,
\varepsilon_4-\varepsilon_3, \varepsilon_5-\varepsilon_1,
\varepsilon_5-\varepsilon_2, \varepsilon_5-\varepsilon_3,
\varepsilon_6-\varepsilon_1, \varepsilon_6-\varepsilon_2,
\varepsilon_6-\varepsilon_3, \varepsilon_7-\varepsilon_6,\right. \cr
&\phantom{J = } \left.\varepsilon_6+\varepsilon_1,
\varepsilon_5+\varepsilon_1,
\varepsilon_5+\varepsilon_2, \varepsilon_4+\varepsilon_1,
\varepsilon_4+\varepsilon_2 \right\}. \cr
}
$$
The placed configuration of boxes corresponding to $(C_0,J_0)$
is as given below.
$$
\matrix{
\beginpicture
\setcoordinatesystem units <0.5cm,0.5cm> point at 6 0       % sets scale
\setplotarea x from 0 to 6, y from -1 to 7    % sets plot size up
\linethickness=0.5pt                          % sets line thickness
%\put{contents of boxes}  at 2.5 -0.5    % 
\put{-2}  at 2.5 1.5    % 
\put{-1}  at 2.5 2.5    % 
\put{-1}  at 3.5 1.5    % 
\put{-1}  at 4.5 0.5    % 
\put{0}  at 2.5 3.5    %  labels contents of boxes
\put{0}  at 3.5 2.5    % 
\put{0}  at 5.5 0.5    % 
\put{1}  at 0.5 6.5    % 
\put{1}  at 1.5 5.5    % 
\put{1}  at 2.5 4.5    % 
\put{2}  at 2.5 5.5    % 
\putrule from 0 7 to 1 7          %            
\putrule from 0 6 to 3 6          %            
\putrule from 1 5 to 3 5          %            
\putrule from 2 4 to 3 4          %            
\putrule from 2 3 to 4 3          %        
\putrule from 2 2 to 4 2          %  draws horizontal lines         
\putrule from 2 1 to 6 1          %           
\putrule from 4 0 to 6 0          %           
\putrule from 0 6 to 0 7        %             
\putrule from 1 5 to 1 7        %             
\putrule from 2 1 to 2 6        %             
\putrule from 3 1 to 3 6        %  draws vertical lines           
\putrule from 4 0 to 4 3        %             
\putrule from 5 0 to 5 1        %  
\putrule from 6 0 to 6 1        %  
\endpicture
&\qquad
&\beginpicture
\setcoordinatesystem units <0.5cm,0.5cm> point at 6 0       % sets scale
\setplotarea x from 0 to 6, y from -1 to 7    % sets plot size up
\linethickness=0.5pt                          % sets line thickness
%\put{numbering of boxes}  at 2.5 -0.5    % 
\put{-7}  at 2.5 1.5    % 
\put{-6}  at 2.5 2.5    % 
\put{-5}  at 3.5 1.5    % 
\put{-4}  at 4.5 0.5    % 
\put{1}  at 2.5 3.5    %  labels contents of boxes
\put{2}  at 3.5 2.5    % 
\put{3}  at 5.5 0.5    % 
\put{4}  at 0.5 6.5    % 
\put{5}  at 1.5 5.5    % 
\put{6}  at 2.5 4.5    % 
\put{7}  at 2.5 5.5    % 
\putrule from 0 7 to 1 7          %            
\putrule from 0 6 to 3 6          %            
\putrule from 1 5 to 3 5          %            
\putrule from 2 4 to 3 4          %            
\putrule from 2 3 to 4 3          %        
\putrule from 2 2 to 4 2          %  draws horizontal lines         
\putrule from 2 1 to 6 1          %           
\putrule from 4 0 to 6 0          %           
\putrule from 0 6 to 0 7        %             
\putrule from 1 5 to 1 7        %             
\putrule from 2 1 to 2 6        %             
\putrule from 3 1 to 3 6        %  draws vertical lines           
\putrule from 4 0 to 4 3        %             
\putrule from 5 0 to 5 1        %  
\putrule from 6 0 to 6 1        %  
\endpicture
&\qquad
&\beginpicture
\setcoordinatesystem units <0.5cm,0.5cm> point at 6 0       % sets scale
\setplotarea x from 0 to 6, y from -1 to 7    % sets plot size up
\linethickness=0.5pt                          % sets line thickness
%\put{standard tableau}  at 2.5 -0.5    % 
\put{4}  at 2.5 1.5    % 
\put{2}  at 2.5 2.5    % 
\put{5}  at 3.5 1.5    % 
\put{6}  at 4.5 0.5    % 
\put{1}  at 2.5 3.5    %  labels contents of boxes
\put{3}  at 3.5 2.5    % 
\put{7}  at 5.5 0.5    % 
\put{-6}  at 0.5 6.5    % 
\put{-5}  at 1.5 5.5    % 
\put{-2}  at 2.5 4.5    % 
\put{-4}  at 2.5 5.5    % 
\putrule from 0 7 to 1 7          %            
\putrule from 0 6 to 3 6          %            
\putrule from 1 5 to 3 5          %            
\putrule from 2 4 to 3 4          %            
\putrule from 2 3 to 4 3          %        
\putrule from 2 2 to 4 2          %  draws horizontal lines         
\putrule from 2 1 to 6 1          %           
\putrule from 4 0 to 6 0          %           
\putrule from 0 6 to 0 7        %             
\putrule from 1 5 to 1 7        %             
\putrule from 2 1 to 2 6        %             
\putrule from 3 1 to 3 6        %  draws vertical lines           
\putrule from 4 0 to 4 3        %             
\putrule from 5 0 to 5 1        %  
\putrule from 6 0 to 6 1        %  
\endpicture
\cr
\hbox{contents of boxes}
&&\hbox{numbering of boxes}
&&\hbox{a standard tableau} \cr
}
$$
\endexample

\subsection{}

Using the above rules one produces a book of
placed configurations corresponding to 
$(\vec\gamma, J)=((\beta_1+C_{\beta_1},J_{\beta_1}),
\ldots,(\beta_r+C_{\beta_r},J_{\beta_r}))$.
A {\it standard tableau} $t$ for this book of configurations is
a filling of the boxes with distinct elements of
$\{-n,\ldots,-1,1,\ldots,n\}$ such that the filling on each page
satisfies the conditions for a standard tableau for that page.
Let
$$\cF^{((C_{\beta_1},J_{\beta_1}),
\ldots,(C_{\beta_r},J_{\beta_r}))}$$
denote the set of such fillings.  

The proof of the following Theorem is similar to
the proof of Theorem (3.5).

\thm  Given a standard tableau $t$ for the book of configurations
$((C_{\beta_1},J_{\beta_1}),
\ldots,(C_{\beta_r},J_{\beta_r}))$ define $w_t\in WC_n$
by $w_t(i)=t({\rm box}_i)$.
Then the map
$$
\matrix{
\cF^{((C_{\beta_1},J_{\beta_1}),
\ldots,(C_{\beta_r},J_{\beta_r}))}
&\longleftrightarrow &\cF^{(\vec\gamma,J)} \cr
t &\longmapsto &w_t \cr }$$
is a bijection.  
\endthm

\vfill\eject
\ 
\bigskip
\centerline{\smallcaps References}
\bigskip

%\medskip
%\item{[Ar]} {\smallcaps S.\ Ariki}, 
%{\it On the decomposition number of the Hecke algebra of $G(m,1,n)$}, 
%FIX THIS REFERENCE.

\medskip
\item{[AL]} {\smallcaps C. Athanasiadis and S. Linusson}, 
{\it A simple bijection for the regions of the Shi arrangement of hyperplanes},
Discrete Math., to appear.

\medskip
\item{[Au]} {\smallcaps M. Aubert}, {\it Dualit\'e dans le groupe de
Grothendieck de la categorie des repr\'esentations lisses de longueur
finie d'un groupe r\'eductif $p$-adique},  Trans. Amer. Soc. {\bf 347}
(1995), 2179--2189.

\medskip
\item{[Bj]} {\smallcaps A. Bj\"orner}, {\it Orderings of Coxeter groups},
{\sl Combinatorics and Algebra (Boulder, Colo. 1983)}, Contemp. Math.
{\bf 34}, Amer. Math. Soc., Providence 1984, 175--195.

\medskip
\item{[Bou]} {\smallcaps N.\ Bourbaki}, 
{\it Groupes et alg\`ebres de Lie, Chapitres 4,5 et 6},
Elements de Math\'ematique, Hermann, Paris 1968.

%\medskip
%\item{[Ca]} {\smallcaps R.W.\ Carter}, 
%{\it Finite Groups of Lie Type: Conjugacy Classes and Complex Characters},
%John Wiley and Sons, 1993.

\medskip
\item{[Fo]} {\smallcaps S. Fomin}, Lecture notes, MIT, Fall 1996.

\medskip
\item{[GK]} {\smallcaps I.M. Gelfand, D. Krob, A.  Lascoux, B. Leclerc, 
V. Retakh, J.-Y. Thibon}, 
{\it Noncommutative symmetric functions}, Adv. Math. {\bf 112}
(1995), 218--348.

\medskip
\item{[He]} {\smallcaps P. Headley}, 
{\it On a family of hyperplane arrangements related to the affine Weyl groups}, 
J. Algebraic Combinatorics {\bf 6} (1997), 331--338. 

%\medskip
%\item{[Ho]} {\smallcaps P.N.\ Hoefsmit}, 
%{\it Representations of Hecke algebras of finite groups with 
%$BN$-pairs of classical type}, 
%Ph.D.\ Thesis, University of British Columbia, 1974. 

%\medskip
%\item{[HO1]} {\smallcaps G.J.\ Heckman and E.M.\ Opdam}, 
%{\it Yang's system of particles and Hecke algebras}, 
%Ann. Math. (2) {\bf 145} (1997), 139--173.

\medskip
\item{[HO]} {\smallcaps G.J.\ Heckman and E.M.\ Opdam}, 
{\it Harmonic analysis for affine Hecke algebras}, 
Current Developments in Mathematics, 1996, International Press, Boston.

\medskip
\item{[Hu]} {\smallcaps J.E.\ Humphreys}, 
{\it Reflection Groups and Coxeter Groups}, 
Cambridge Studies in Mathematics {\bf 29}, Cambridge Univ. Press 1990.

\medskip
\item{[K1]} {\smallcaps S. Kato}, 
{\it Irreducibility of principal series representations for Hecke 
algebras of affine type},  
J. Fac. Sci. Univ. Tokyo sec 1A {\bf 28} (1981), 929--943.

\medskip
\item{[K2]}  {\smallcaps S. Kato}, {\it Duality for representations of
a Hecke algebra}, Proc. Amer. Math. Soc. {\bf 119} (1993), 941--946.

%\medskip
%\item{[KL]} {\smallcaps D.\ Kazhdan and G.\ Lusztig}, 
%{\it Proof of the Deligne-Langlands conjecture for Hecke algebras}, 
%Invent. Math. {\bf 87} (1987), 153--215.

\medskip
\item{[KZ]} {\smallcaps H. Knight and A.V. Zelevinsky}, {\it Representations
of quivers of type A and the multisegment duality}, Adv. in Math., to appear.

\medskip
\item{[La1]} {\smallcaps A. Lascoux, B. Leclerc, J.-Y. Thibon}, 
{\it Ribbon tableaux,
Hall-Littlewood functions, quantum affine algebras, and unipotent varieties}, 
J. Math. Phys. {\bf 38} (1997), 1041--1068.

\medskip
\item{[La2]} {\smallcaps A. Lascoux, P. Pragacz}, {\it Ribbon Schur functions}, 
European J. Combin. {\bf 9} (1988), 561--574.

\medskip
\item{[Li1]} {\smallcaps P. Littelmann}, {\it Paths and root operators in
representation theory}, Ann. of Math. (2) {\bf 142} (1995), 499--525.

\medskip
\item{[Li2]} {\smallcaps P. Littelmann}, {\it The path model for
representations of symmetrizable Kac-Moody algebras}, Proceedings of the
International Congress of Mathematicians, Vol. 1 (Z\"uŸrich, 1994), 298--308,
BirkhaŠuser, Basel, 1995. 

\medskip
\item{[Mac]} {\smallcaps I.G. Macdonald}, 
Symmetric functions and Hall polynomials, Second edition, 
Oxford Mathematical Monographs, Oxford Univ. Press, New York, 1995.

%\medskip
%\item{[Mac2]} {\smallcaps I.G.\ Macdonald}, 
%{\it Affine Hecke algebras and orthogonal polynomials}, 
%S\'eminaire Bourbaki, 47\`eme ann\'ee, ${\rm n}^{\rm o}$ 797, 1994--95.

\medskip
\item{[Mat]} {\smallcaps H. Matsumoto}, 
{\it Analyse harmonique dans les systems de Tits bornologiques de type affine},
Lect. Notes in Math. {\bf 590}, Springer-Verlag, 1977.

\medskip
\item{[MW]} {\smallcaps Moeglin and Waldspurger},
{\it L'involution de Zelevinski}, J. Reine Angew. Math. {\bf 372} (1986),
136--177. 

\medskip
\item{[Ra1]} {\smallcaps A. \ Ram}, 
{\it Calibrated representations of affine Hecke algebras}, in preparation.

\medskip
\item{[Ra2]} {\smallcaps A.\ Ram}, 
{\it Irreducible representations
of rank two affine Hecke algebras}, in preparation.

\medskip
\item{[Ra3]} {\smallcaps A.\ Ram}, 
{\it Skew shape representations are irreducible}, in preparation.

\medskip
\item{[RR1]} {\smallcaps A.\ Ram and J. Ramagge}, 
{\it  Jucys-Murphy elements come from affine Hecke algebras}, in preparation.

\medskip
\item{[RR2]} {\smallcaps A.\ Ram and J. Ramagge}, {\it Calibrated
representations and the $q$-Springer correspondence}, in preparation.

\medskip
\item{[Ro]} {\smallcaps F. Rodier}, 
{\it Decomposition de la s\'erie principale des groups reductifs $p$-adic}, 
in Noncommutative Harmonic analysis and
Lie Groups, Lect. Notes in Math. {\bf 880}, Springer-Verlag, 1981.

\medskip
\item{[Rg]} {\smallcaps J. Rogawski}, {\it On modules over the Hecke algebra
of a $p$-adic group}, Invent. Math. {\bf 79} (1985), 443--465.

\medskip
\item{[Sh1]} {\smallcaps J.-Y. Shi}, {\it The number of $\oplus$-sign types},
Quart. J. Math. Oxford Ser. (2) {\bf 48} (1997), 93--105.

\medskip
\item{[Sh2]} {\smallcaps J.-Y. Shi}, {\it  Left cells in affine Weyl
groups}, T\^ohoku Math. J. (2) {\bf 46} (1994), no. 1, 105--124. 

\medskip
\item{[Sh3]} {\smallcaps J.-Y. Shi}, 
{\it  Sign types corresponding to an affine Weyl group}, 
J. London Math. Soc. (2) {\bf 35} (1987), 56--74.

\medskip
\item{[ST]} {\smallcaps L. Solomon and H. Terao}, 
{\it The double Coxeter arrangement}, preprint 1997.
 
\medskip
\item{[St1]} {\smallcaps R. Stanley}, 
{\it Hyperplane arrangements, interval orders, and trees}, 
Proc. Nat. Acad. Sci. U.S.A. {\bf 93} (1996), 2620--2625.

\medskip
\item{[St2]} {\smallcaps R. Stanley}, 
{\it Hyperplane arrangements, parking functions and tree inversions}, 
in {\it Mathematical Essays in Honor of Gian-Carlo Rota}, 
Birkh\"auser, Boston, 1998.

%\medskip
%\item{[Sng]} {\smallcaps R.\ Steinberg}, 
%{\it Endomorphisms of linear algebraic groups}, 
%Memoirs of the American Mathematical Society, {\bf 80} (1968) 1--108.

\medskip
\item{[Xi]} {\smallcaps N. Xi}, 
{\it  Representations of affine Hecke algebras}, 
Lect. Notes in Math. {\bf 1587}, Springer-Verlag, Berlin, 1994.

\medskip
\item{[Wz]} {\smallcaps H. Wenzl}, {\it Hecke algebras of type $A_n$ and
subfactors}, Invent. Math. {\bf 92} (1988), 349--383.

\medskip
\item{[Y]} {\smallcaps A. Young}, {\it On quantitative substitutional analysis
(sixth paper)}, Proc. London. Math. Soc. (2) {\bf 34} (1931), 196--230.

\medskip
\item{[Ze]} {\smallcaps A. Zelevinsky}, 
{\it Induced Representations of $\gp$-adic groups II: 
On Irreducible Representations of $GL(n)$}, 
Ann. Scient.  \'Ec. Norm.  Sup. $4^{\rm e}$ s\`erie, {\bf 13} (1980), 165--210. 

\vfill\eject
\end

\subsection{Decomposition of the principal series module}

Suppose $M(\rho)$ is the
principal series module corresponding to $\rho$. As an
$H$-module, $M(\rho)$ is isomorphic to the left regular representation. In
particular, $M(\rho)$ has basis $\{T_wv_\rho\}_{w\in W}$. In the Grothendieck ring
we have 
$$
\left[ M(\rho)\right] = \sum_{J\subseteq P(\rho)} \left[ H^{(\rho,J)}\right].$$
This identity is a module theoretic version of (?) for the case when
$\gamma=\rho$. 

In type $A$ the identity (??) has appeared in the form
$$
\CC\left[ S_n\right]\cong \bigoplus_{\theta} S^{\theta}
$$ 
where the sum is over all ribbons $\theta$ of length $n$ and
$S^{\theta}$ denote the corresponding skew shape representations of
$S_n$. 

Finally, the identity has also appeared in the symmetric function literature as 
$$
h_1^n = \sum_{\theta} s_{\theta}
$$
where the sum is over all ribbons $\theta$ of length $n$, $h_1=x_1+\cdots+x_n$
is the first complete symmetric function and
$s_{\theta}$ are skew Schur functions [Mac I \S5 Ex.~21b].

\subsection Row reading talbeaux

Let
Suppose $j>i$ and 
$\varepsilon_j-\varepsilon_i\in\overline{J}$. Then 
$\varepsilon_j-\varepsilon_i=
(\varepsilon_j-\varepsilon_{j_1})+(\varepsilon_{j_1}-\varepsilon_{j_2})
+\cdots+(\varepsilon_{j_r}-\varepsilon_i)$ where each of the
summands is in $J$.  This means that ${\rm box}_j$ is northwest of
${\rm box}_{j_1}$ is northwest of ${\rm box}_{j_2}$ is northwest of
$\ldots$ is notrthwest of ${\rm box}_i$.  Thus the 
position of ${\rm box}_j$ relative to ${\rm box}_i$ is somewhere in the shaded
region of 
$$
\beginpicture
\setcoordinatesystem units <0.75cm,0.75cm>         % sets scale
\setplotarea x from 0 to 4, y from 0 to 5    % sets plot size up
\linethickness=0.5pt                          % sets line thickness
\put{\vdots} at 0 4.5 
\put{\vdots} at 4 4.5 
\put{${\scriptstyle{\rm box}_i}$} at 3.5 0.5 
\putrule from 1 3 to 2 3          %  draws horizontal lines         
\putrule from 2 2 to 3 2          %           
\putrule from 3 1 to 4 1          %           
\putrule from 3 0 to 4 0          %           
\putrule from 0 4 to 1 4        %             
\putrule from 1 3 to 1 4        %             
\putrule from 2 2 to 2 3        %  draws vertical lines           
\putrule from 3 0 to 3 2        %             
\putrule from 4 0 to 4 4        %  
\vshade 0 4 5   1 4 5 / 
\vshade 1 3 5   2 3 5 /
\vshade 2 2 5   3 2 5 /
\vshade 3 1 5   4 1 5 /           
\endpicture
$$
and
$\varepsilon_j-\varepsilon_i\in\overline{J}^c$ if ${\rm box}_j$ is in the
following shaded region.
$$
\beginpicture
\setcoordinatesystem units <0.75cm,0.75cm>         % sets scale
\setplotarea x from 0 to 5, y from 0 to 6    % sets plot size up
\linethickness=0.5pt                          % sets line thickness
\put{\vdots} at 1 5.5 
\put{\vdots} at 3 0.5 
\put{${\scriptstyle{\rm box}_i}$} at 0.5 3.5 
\putrule from 0 4 to 1 4          %            
\putrule from 0 3 to 1 3          %  draws horizontal lines         
\putrule from 1 2 to 2 2          %           
\putrule from 2 1 to 3 1          %           
\putrule from 0 3 to 0 4        %             
\putrule from 1 2 to 1 5        %             
\putrule from 2 1 to 2 2        %  draws vertical lines           
\vshade 1 2 6   2 2 6  /
\vshade 2 1 6   3 1 6  /
\vshade 3 0 6   4 0 6 /           
\endpicture
$$

For the purposes of contradiction assume that
$j>i$, $l>k$, with 
$\varepsilon_j-\varepsilon_i,\varepsilon_l-\varepsilon_k\in\overline{J}^c$ and
$(\varepsilon_j-\varepsilon_i)+(\varepsilon_l-\varepsilon_k)
\in R^+\setminus\overline{J}^c$. 
In particular, either $j=k$ or $i=l$. 
Suppose $i=l$. Hence
$\varepsilon_j-\varepsilon_i,\varepsilon_i-\varepsilon_k\in\overline{J}^c$ and
$\varepsilon_j-\varepsilon_k\not\in\overline{J}^c$. 
In other words,
$\varepsilon_j-\varepsilon_i,\varepsilon_i-\varepsilon_k\in\overline{J}^c$ and
$\varepsilon_j-\varepsilon_k\in\overline{J}$. 
By considering their relative
positions in the shape we find this to be impossible. 
Similarly we cannot have $j=k$. 
Thus $\overline{J}^c$ is closed and $\cF^{(\gamma,J)}$ has a unique minimal
element.